\documentclass[12pt,a4paper,twoside]{amsart}

\usepackage[T1]{fontenc}
\usepackage{lmodern}
\usepackage{amsmath,amssymb,amsthm,mathrsfs}
\IfFileExists{stmaryrd.sty}{\usepackage{stmaryrd}}{%
\newcommand{\llbracket}{\mathopen{[\![}}%
\newcommand{\rrbracket}{\mathclose{]\!]}}%
}
\usepackage[margin=1in]{geometry}
\usepackage{hyperref}
\newtheorem{theorem}{Theorem}[section]

\newtheorem{lemma}[theorem]{Lemma}
\newtheorem{claim}[theorem]{Claim}
\newtheorem{corollary}[theorem]{Corollary}

\theoremstyle{definition}
\newtheorem{definition}[theorem]{Definition}
\newtheorem{hypothesis}[theorem]{Hypothesis}

\theoremstyle{remark}

\newtheorem{notation}[theorem]{Notation}

\DeclareMathOperator{\Col}{Col}
\DeclareMathOperator{\cf}{cf}
\DeclareMathOperator{\dom}{dom}
\DeclareMathOperator{\ran}{ran}
\DeclareMathOperator{\lh}{lh}
\DeclareMathOperator{\Ult}{Ult}

\DeclareMathOperator{\Lev}{Lev}
\DeclareMathOperator{\RO}{RO}
\DeclareMathOperator{\Succ}{Succ}
\DeclareMathOperator{\stem}{stem}
\DeclareMathOperator{\supp}{supp}
\DeclareMathOperator{\tail}{tail}

\newcommand{\Q}{\mathbb Q}
\newcommand{\Pp}{\mathbb P}
\newcommand{\R}{\mathbb R}

\newcommand{\restr}{\restriction}
\newcommand{\forces}{\Vdash}
\newcommand{\Ord}{\operatorname{Ord}}

\newcommand{\concat}{\mathbin{{}^\frown}}
\newcommand{\pow}{\mathcal P}
\newcommand{\cT}{\mathcal T}

\newcommand{\bb}[1]{\llbracket #1\rrbracket}

\newcommand{\Sforce}{\mathbb S}
\newcommand{\res}{\restr}

\newcommand{\PP}{\mathbb P_{\mathrm{no}}}

\numberwithin{equation}{section}

\begin{document}
\title[Adding minimal fresh subsets to  $\aleph_\omega$]
      {Adding minimal fresh subsets to  $\aleph_\omega$  }

\author[M. Gitik]{Moti Gitik}
\address{School of Mathematical Sciences, Tel Aviv University, Tel Aviv,
Israel.}
\email{gitik@post.tau.ac.il}

\author[M.  Golshani]{Mohammad Golshani}
	\address{Mohammad Golshani, School of Mathematics, Institute for Research in Fundamental Sciences (IPM), P.O.\ Box:
		19395--5746, Tehran, Iran.}
	\email{golshani.m@gmail.com}

\thanks{The second author's research has been supported by a grant from IPM (No. 1405030417).}


\begin{abstract}
Assuming the existence of a measurable cardinal, we produce
generic extensions $V_1\subsetneq V_2$ of the universe with
the same cardinals, cofinalities, and bounded subsets of $\aleph_\omega$, such
that $V_2$ is generated over
$V_1$ by a  cofinal $\omega$-sequence in $\aleph_\omega$, and $V_2$ is a
minimal generic extension over $V_1$, i.e., every set of ordinals in
$V_2\setminus V_1$ generates $V_2$ over $V_1$.
\end{abstract}
\maketitle
\setcounter{section}{-1}

\section{Introduction}

Suppose $M\subseteq N$ are transitive models of ZFC with the same
ordinals, and $\kappa$ is a cardinal of $N$. A set $A\subseteq\kappa$
is called \emph{fresh over $M$} if $A\notin M$ and
$A\cap\alpha\in M$ for all $\alpha<\kappa$. We say that $N$ is
\emph{minimal over $M$} if, whenever $A\in N\setminus M$ is a set of
ordinals, one has $N=M[A]$.

It is well known that Prikry forcing yields many intermediate submodels, see \cite{gitikkanovei} and \cite{Sureson}. Later work by Koepke, R\H{a}sch, and Schlicht \cite{KoepkeRaeschSchlicht}, as well as Yang \cite{Yang}, produced variants of Prikry forcing that yield minimal extensions of the universe. In these models, all cardinals are preserved, and hence the minimality is with respect to a large cardinal. In this paper, we show that one can   pushed down the result to small accessible cardinals.
We fist prove the following result.

\begin{theorem}\label{thm:main}
Assume
$\langle\kappa_n:n<\omega\rangle$ is an increasing sequence of
measurable cardinals such that for every $n>0$, $\kappa_n$ carries at
least $\kappa_{n-1}$ distinct normal measures, and let
$\kappa=\sup_{n<\omega}\kappa_n$. Then there are generic extensions
$V_1\subsetneq V_2$ of $V$ satisfying the following conditions:
\begin{enumerate}
\item $V_1$ and $V_2$ have the same cardinals and cofinalities, and
$\pow(\alpha)^{V_1}=\pow(\alpha)^{V_2}$ for every $\alpha<\kappa$;
\item $\kappa=\aleph_\omega^{V_1}=\aleph_\omega^{V_2}$;
\item there is an increasing cofinal sequence
$C=\langle\xi_n:n<\omega\rangle$ in $\kappa$ such that $V_2=V_1[C]$;
\item for every set of ordinals $X\in V_2\setminus V_1$,
$V_1[X]=V_2$.
\end{enumerate}
\end{theorem}
Note that the hypotheses of the theorem follow from the assumption of the existence of an increasing sequence
$\langle\kappa_n:n<\omega\rangle$ of measurable cardinals.

The proof of the theorem is in essentially two steps.  First,  from the above
sequence of measures, we collapse enough cardinals between  $\kappa_n$'s to make $\kappa=\aleph_\omega$ in such a way that the
measures generate normal ideals whose positive-set orders are highly strategically closed.   Second, we modify the forcing construction
 from Yang Sen \cite[Section 4]{Yang} and define tree Prikry forcing notions built  from these ideals
to add a fresh cofinal branch that yields a minimal extension.
Prikry type forcing notions using suitable ideals were already discussed, see for example \cite{{BNI}}
and \cite{Sakai}.

We then show that one can prove the same result, using just one measurable cardinal.
\begin{theorem}\label{thm:main2}
Assume
$\kappa$ is a
measurable cardinal. Then there are generic extensions
$V_1\subsetneq V_2$ of $V$ satisfying the following conditions:
\begin{enumerate}
\item $V_1$ and $V_2$ have the same cardinals and cofinalities, and
$\pow(\alpha)^{V_1}=\pow(\alpha)^{V_2}$ for every $\alpha<\kappa$;
\item $\kappa=\aleph_\omega^{V_1}=\aleph_\omega^{V_2}$;
\item there is an increasing cofinal sequence
$\vec\eta=\langle\xi_n:n<\omega\rangle$ in $\kappa$ such that $V_2=V_1[\vec\eta]$;
\item for every set of ordinals $X\in V_2\setminus V_1$,
$V_1[X]=V_2$.
\end{enumerate}
\end{theorem}
The proof of this theorem uses a variant of Prikry forcing which on the one hand, adds an $\omega$-Prikry sequence with collapses between them to make $\kappa=\aleph_\omega$, and on the other hand, defines a variant of \cite{KoepkeRaeschSchlicht} to add a minimal $\omega$-cofinal sequence through the resulting model.

\section{Forcing a minimal  subset of \texorpdfstring{$\aleph_\omega$}{aleph omega} using $\omega$-many measurable cardinals}
In this section we prove Theorem \ref{thm:main}. Starting from the hypotheses of the theorem, we first introduce a
system of normal ideas which has several nice properties, see Theorem \ref{thm:preparation}, and then use this system of ideals to define a variant of Prikry forcing which gives a minimal generic extension.
\subsection{A nice system of Normal ideals}
In this subsection starting from hypotheses of the main theorem, we find a generic extension which makes $\kappa=\aleph_\omega$
and in the extension there is a nice system of normal measures. Similar ideas were developed in \cite{BNI}.

Assume  $\rho$ is a regular uncountable cardinal. Recall that a
 \emph{normal ideal} on $\rho$ is a proper ideal $I\subseteq\pow(\rho)$
which is $\rho$-complete,
and is closed under diagonal union, i.e., if $\langle A_\xi:\xi<\rho\rangle$ is a sequence from $I$,
then their diagonal union
\[
 \nabla_{\xi<\rho}A_\xi
 :=\{\alpha<\rho:(\exists\xi<\alpha)\ \alpha\in A_\xi\}
\]
is also in $I$.
As usual, we let $I^+=\pow(\rho)\setminus I$ and
$I^*=\{A\subseteq\rho:\rho\setminus A\in I\}$.
The following is well-known.

\begin{lemma}\label{lem:normalfacts}
Assume $I$ is a normal ideal on $\rho$.
\begin{enumerate}
\item every bounded subset of $\rho$ is in $I$, hence, every member of $I^+$
has cardinality $\rho$.
\item If $A\in I^+$ and $A=\bigcup_{i<\eta}A_i$, where $\eta<\rho$,
then $A_i\in I^+$ for some $i<\eta$.
\item If $A\in I^+$ and $f:A\setminus\{0\}\to\rho$ is regressive,
then $f^{-1}\{\xi\}\in I^+$ for some $\xi<\rho$.
\end{enumerate}
\end{lemma}
\begin{definition}\label{def:goodideal}
Assume $\lambda\leq\rho$ are  regular cardinals and $\rho$ is
uncountable. A normal ideal $I$ on $\rho$ is
\emph{$(\lambda,\rho)$-good} if the forcing notion
$(I^+,\subseteq)$ is $\lambda$-strategically closed. When $\rho$ is
clear from context, we also say that $I$ is $\lambda$-good.
\end{definition}

We now state and prove some general results that show that one can produce good ideas in forcing extensions.

\begin{lemma}\label{lem:generatedideal}
Assume $U$ is a normal measure on $\rho$  and $\Pp$ is
$\rho$-c.c. If $G$ is $\Pp$-generic over $V$, then
\[
 I_U^G=\{A\subseteq\rho:(\exists B\in U)\ A\cap B=\varnothing\}
\]
is a normal ideal on $\rho$ in $V[G]$.
\end{lemma}
\begin{proof}
It is clear that $I_U^G$ is an ideal on $\rho.$
We first show that it is  $\rho$-complete.  Assume that  $\delta\leq\rho$ and
$\langle A_\xi:\xi<\delta\rangle$ is a sequence in $V[G]$ of members
of $I_U^H$.
Let $p \in G$ and the sequence $\langle \dot B_\xi: \xi < \delta \rangle$ of $\Pp$-names be such that
$p$ forces ``for all $\xi < \delta,$ $\dot B_\xi \in U$ and  $\dot A_\xi\cap\dot B_\xi=\varnothing$''.

For each $\xi$, let  $\mathcal A_\xi$ and $\langle B_{\xi,q}: q\in\mathcal A_\xi   \rangle$ be such that:
\begin{itemize}
\item  $\mathcal A_\xi$  is a
 maximal antichain below $p$,

 \item For each $q\in\mathcal A_\xi$, $B_{\xi,q}\in U$, and    $q$ forces ``$\dot B_\xi=B_{\xi,q}\in U$''.
\end{itemize}
Note that by our assumption, $|\mathcal A_\xi| < \rho$, and hence
$
 C_\xi=\bigcap_{q\in\mathcal A_\xi}B_{\xi,q}\in U.
$
Clearly
$p\forces\check C_\xi\subseteq\dot B_\xi$.
Now $C=\bigcap_{\xi<\delta}C_\xi \in U$
and it is disjoint from $\bigcup_{\xi<\delta}A_\xi$.

To show $I_U^G$ is closed under diagonal union, let $\langle A_\xi: \xi < \rho  \rangle$
be  a sequence in $V[G]$ from $I_U^G$. Proceed exactly as above and produce the sets $C_\xi \in U$ for $\xi < \rho$.
Then the diagonal intersection
$
 C=\{\alpha<\rho:(\forall\xi<\alpha)\ \alpha\in C_\xi\}
$
belongs to $U$ and is disjoint from
$\nabla_{\xi<\rho}A_\xi$.
\end{proof}
The next lemma is well-know.
\begin{lemma}\label{lem:smallforcing}
Assume $U$ is a normal measure on $\rho$  and $|\Pp|< \rho$. Let $G$ be $\Pp$-generic over $V$.
Then in $V[G]$, the filter generated by $U$  is
a normal measure on $\rho$.
Furthermore, distinct normal measures have distinct generated extensions.
\end{lemma}

\begin{lemma}\label{lem:collapsedideal}
Assume $U$ is a normal measure on $\rho$  and  $\mu<\rho$ is regular
and uncountable. Let $G$ be
$\Pp=\Col(\mu,<\rho)$-generic over $V$. Then
in $V[G]$, the ideal
\[
 I=\{A\subseteq\rho:(\exists B\in U)\ A\cap B=\varnothing\}
\]
is $(\mu, \rho)$-good
\end{lemma}
\begin{proof}
The Levy collapse  $\Pp=\Col(\mu,<\rho)$ is $\rho$-c.c., hence by
Lemmas~\ref{lem:generatedideal}, $I$ is normal. We show that $(I^+, \subseteq)$ is $\mu$-strategically closed.

Let $j:V\to M=\Ult(V,U)$ be the ultrapower embedding  and set
$
 \Q=\Col(\mu,[\rho,j(\rho)))^M.
$
Clearly  $j(\Pp)\simeq \Pp\times \Q$ and   $\Q$ is $<\mu$-closed in $v[G]$.
Let $H$ be $\Q$-generic over $V[G]$ and define $j_{G, H}: V[G] \longrightarrow M[G\times H]$
by
\[
 j_{G, H}(\dot x[G])=j(\dot x)[{G\times H}].
\]
Then by Silver's theorem, $j_{G, H}$ is a well-defined elementary embedding extending $j$, and
$j_{G, H}(G)=G\times H$.
Let $\mathbb B=\RO(\Q)^{V[G]}$, and
for $A=\dot A[G] \subseteq\rho$, set
\[
 e(A)=\bigl\|\check\rho\in
          j(\dot A)[{G\times \dot H}]\bigr\|_{\mathbb B}.
\]
Let also $[r]_{\mathbb B}$ denote the non-zero Boolean value of a condition $r\in \Q$,
namely $[r]_{\mathbb B}=\bigl\|   r\in \dot H  \bigr\|_{\mathbb B}$.
\begin{claim}\label{claim:boolproperties}
The map $e$ is well-defined, is a Boolean homomorphism, and for every $\delta<\rho$
 and $\langle A_\xi:\xi<\delta\rangle \in V[G]$, we have
$
 e\left(\bigcap_{\xi<\delta}A_\xi\right)
       =\bigwedge_{\xi<\delta}e(A_\xi).
$
Moreover, the kernel of $e$ is $I$.
\end{claim}
\begin{proof}
It is easy to see that $e$ is a well-defined  Boolean homomorphism
and it satisfies $
 e\left(\bigcap_{\xi<\delta}A_\xi\right)
       =\bigwedge_{\xi<\delta}e(A_\xi),
$  for every
 $\delta<\rho$
 and $\langle A_\xi:\xi<\delta\rangle \in V[G]$.

We now show that the kernel of $e$ is  $I$.
first suppose that $A\in I$. Let $B\in U$ be disjoint from $A$.
Then  $\rho\in j(B)$, so $\rho\notin j_{G, H}(A)$, and hence
 $e(A)=0$.
Conversely, let $A=\dot A[G]$ and assume that $e(A)=0$. This means that
$\bigl\|\check\rho\in j(\dot A)[{G\times H}]\bigr\|_{\mathbb B}=0$.
Using the forcing theorem and homogeneity of the forcing, let us assume that
$
 (1_{\Pp},1_\Q)\forces^M_{\Pp\times \Q}\check\rho\notin j(\dot A).
$
By elementarity, the set
\[
 B=\{\alpha<\rho:  1_\Pp \forces_\Pp\check\alpha\notin\dot A\}
\]
belongs to $U$, and it is clearly disjoint from $A$.
Thus $A\in I$.
\end{proof}

\begin{claim}\label{claim:represent}
For every $r\in \Q$ there exists a set $E_r\subseteq\rho$ in $v[G]$ such that
$e(E_r)=[r]_{\mathbb B}$.
\end{claim}
\begin{proof}
Assume $r\in \Q$. Let $f:\rho\to \Pp$ be such that $j(f)(\rho)=(1_\Pp, r)$. Set
$
 E_r=\{\alpha<\rho:f(\alpha)\in G\}.
$
Then we have
\[
 \rho\in j_{G, H}(E_r)
\Longleftrightarrow j(f)(\rho)\in G\times H
\Longleftrightarrow r\in H.
\]
Taking Boolean values, we get
$e(E_r)=\bigl\|   \rho\in j_{G, H}(E_r)   \bigr\|= \bigl\|   r\in \dot H  \bigr\|=[r]_{\mathbb B}$
as required.
\end{proof}
We now show that $(I^+, \subseteq)$ is $\mu$-strategically closed.
Fix $A\in I^+$ and $\eta<\mu$.
We define a winning strategy for player II for a game of length $\eta$ as follows.
At step $0$, choose $r_0\in \Q$ with
$0<[r_0]_{\mathbb B}\leq e(A)$ and play $C_0=A\cap E_{r_0}$, where $E_{r_0}$ is as the above claim.
Then $e(C_0)=[r_0]_{\mathbb B}>0$, hence $C_0 \in I^+$.
At step $\xi$ player I chooses some $A_\xi \in I^+$. Player II first  chooses an auxiliary
$r_\xi\in \Q$  such that
$ [r_\xi]_{\mathbb B}\leq e(A_\xi),$ and
 $ r_\xi\leq r_\zeta$  for all $\zeta<\xi.$
Then he  plays with
$
 C_\xi=A\cap\bigcap_{\zeta<\xi}A_\zeta\cap E_{r_\xi}.
$
Since
$e(C_\xi)=[r_\xi]_{\mathbb B}>0$,
we have $C_\xi\in I^+$.
This defines a winning strategy for player II.
\end{proof}

Assume now that
$\langle\kappa_n:n<\omega\rangle$ is an increasing sequence of measurable cardinals such that,
for every $n>0$, $\kappa_n$ carries  $\kappa_{n-1}$ distinct
normal measures. Set
$
 \kappa=\sup_{n<\omega}\kappa_n.
$
For each $n>0$ choose distinct normal measure $\langle U_{n, \alpha}: \alpha < \kappa_{n-1}          \rangle$
on $\kappa_n$ and let $\langle X_{n, \alpha}: \alpha < \kappa_{n-1}        \rangle$ be a sequence of pairwise disjoint sets such that
$X_{n, \alpha} \in U_{n, \alpha}$ and $X_{n, \alpha} \subseteq (\kappa_{n-1}, \kappa_n)$.
Let also $U_{0}$ be a normal measure  on $\kappa_0$.
Let
\[
\Pp_\omega= \langle \Pp_n,\dot \Q_n:n<\omega\rangle
\]
be the full support iteration  such that  $\Pp_0$ is the trivial forcing notion, and
$1_{\Pp_n}\forces$``$\Q_n=\dot{\Col}(\lambda_n,<\kappa_n)$'', where
\[
 \lambda_n=
 \begin{cases}
 \omega_1,&n=0,\\
 \kappa_{n-1}^+,&n>0.
 \end{cases}
\]
Let $G$ be $\Pp_\omega$-generic over $V$.
Then the following hold in $V[G]:$
\[
 \kappa_n=\aleph_{2n+2},\qquad
 \kappa=\aleph_\omega,\qquad
 \lambda_0=\omega_1,\qquad
 \lambda_n=\kappa_{n-1}^+\quad(n>0).
\]
For  each $n<\omega$ let $G_n$ be the resulting $\Pp_n$-generic filter over $V$.
Then in  $V[G_{n+1}],$ we have $\kappa_i=\aleph_{2i+2}$ for every $i\leq n$, and
the above lemmas show that in $V[G_{n+1}],$ each selected measure $U_{n, \alpha}$ gives a
$(\lambda_n, \kappa_n)$-good ideal $I_{n, \alpha}$ such that $\kappa_n \setminus X_{n, \alpha} \in I_{n, \alpha}$ As the rest of the forcing adds no new subsets to $\kappa_n$,
each $I_{n, \alpha}$ remains a $(\lambda_n, \kappa_n)$-good ideal in $V[G]$. Similarly we get a
$(\lambda_0, \kappa_0)$-good ideal $I_{0}$  on $\kappa_0$, derived from $U_0$.
Putting all the above results together, we get the following theorem.

\begin{theorem}\label{thm:preparation}
Assume $\text{GCH}$ holds and
$\langle\kappa_n:n<\omega\rangle$ is an increasing sequence of
measurable cardinals such that for every $n>0$, $\kappa_n$ carries at
least $\kappa_{n-1}$ distinct normal measures, and let
$\kappa=\sup_{n<\omega}\kappa_n$.  Then there is a set-generic
extension $W$ of $V$ such that:
\begin{enumerate}
\item $\kappa=\aleph_\omega,$

\item $\kappa_n=\aleph_{2n+2},$
 $\lambda_0=\omega_1,$
$ \lambda_n=\kappa_{n-1}^+$, for each $n>0$.

\item There is a $(\lambda_0,\kappa_0)$-good ideal $I_{0,*}$ on $\kappa_0$.

\item For $n>0$ and $\alpha<\kappa_{n-1}$ there is a
$(\lambda_n,\kappa_n)$-good ideal $I_{n,\alpha}$ on $\kappa_n$.

\item For $n>0$ there are pairwise disjoint sets
$\langle X_{n,\alpha}:\alpha<\kappa_{n-1}\rangle$ such that
$
 X_{n,\alpha}\subseteq(\kappa_{n-1},\kappa_n),$ and
$ X_{n,\alpha}\in I_{n,\alpha}^*,$

\item There is a $(\lambda_0,\kappa_0)$-good ideal $I_{0}$ on $\kappa_0$.
\end{enumerate}
\end{theorem}

\subsection{The ideal-tree forcing and its basic properties}
We now introduce our main forcing construction and present some of its properties.
From now on, we assume that the hypotheses (1)--(6)  of Theorem \ref{thm:preparation} hold. The forcing we define is motivated by \cite{Yang}, where normal measures used there are replaced by the ideals found in Theorem \ref{thm:preparation}.
Let $\cT$ be the tree whose nodes are finite sequences $t$ such that
\begin{itemize}
\item $t(i)<\kappa_i,$ for  $i<\lh(t)$,

\item $ t(i)\in X_{i,t(i-1)}$, for every $0<i<\lh(t)$.
\end{itemize}
For $t\in\cT$, we define $I_t$ as
\[
 I(t)=
 \begin{cases}
 I_{0},&\lh(t)=0,\\
 I_{n,t(n-1)},&\lh(t)=n>0.
 \end{cases}
\]
We now define our main forcing notion.
\begin{definition}\label{def:forcing}

A condition in $\Pp$ is a pair $p=(s^p, T^p)$, where $T^p$ is a downward-closed subtree  of $\cT$,
$s^p=\stem(T^p)$ is the stem of $T^p$ such that every node of $T^p$ is comparable with $s^p$ and for every
$s^p\subseteq t\in T^p$,
\[
 \Succ_{T^p}(t):=\{\alpha:t\concat\langle\alpha\rangle\in T^p\}
 \in I(t)^+.
\]
For $p, q\in\Pp$, we say $q$ is an extension of $p$, and write $q \leq p$ when $T^q\subseteq T^p$.  We also define the Prikry order $\leq^*$ by
$q\leq^* p$ iff $q\leq p\ \text{ and }\ s^q=s^p.$
\end{definition}
Note that $1_{\Pp}=(\emptyset, \cT)$ is the largest element of $\Pp.$
\begin{notation}
Assume $p=(s, T) \in \Pp$.
\begin{enumerate}
\item For  $s\subseteq t\in T$, set
\begin{enumerate}
\item $ T\restr t=\{u\in T:u\subseteq t\text{ or }t\subseteq u\},$

\item $ \Lev_k(T)=\{t\in T:\lh(t)=k\}.$
\end{enumerate}
\item Let $[T]$ denote the set of branches through $T$. Thus
 $[T]=\{b\in\prod_{n<\omega}\kappa_n:
(\forall k<\omega)\ b\restr k\in T\}$.

\item If $s\subseteq t\in T$, then $p[t]=(t, T \restr t)$.
\end{enumerate}

\end{notation}
\begin{definition}\label{def:front}
Assume $p=(s, T) \in \Pp$. We say $F\subseteq T$ is a \emph{front} of $p$ if it satisfies the following conditions:
\begin{itemize}
\item Each node in $F$ extends $s$,
\item elements of $F$ are pairwise incompatible  under proper extension,
\item $F$ meets every branch of $T$.
\end{itemize}
For fronts $F,F'$ of $p$, we say $F\prec F'$ if every member of $F'$
properly extends a member of $F$.
\end{definition}
We have the following easy amalgamation lemma.
\begin{lemma}\label{lem:gluing}
Assume $p=(s, T) \in \Pp$, $F$ is a front of $p$,  $T_t\leq^*T\restr t$ for each
$t\in F$ and
$ S=\bigcup_{t\in F}T_t.$ Then
\begin{enumerate}
\item $q=(s, S)$
is a condition, $q\leq^*p$
and $S\restr t=T_t$ for $t\in F$.
\item If $u\in T$ is a proper initial segment of a member of $F$, then
$\Succ_S(u)=\Succ_T(u)$.

\item The set $\{(u, S \restr u): \exists t \in F,~u \subseteq t               \}$is pre-dense below $p$.
\end{enumerate}
\end{lemma}
The next finite homogeneity lemma follows from Lemma~\ref{lem:normalfacts}(ii).
\begin{lemma}\label{lem:homogeneous}
Assume $p=(s, T)\in\Pp$, $m=\lh(\stem(T))$, and $m\leq\ell<\omega$.
For every coloring $c:\Lev_\ell(T)\to r$, where $0<r<\omega$, there
is $S\leq^*T$ such that $c\restr\Lev_\ell(S)$ is constant.
\end{lemma}
\begin{lemma}[Strong Prikry property]\label{lem:strong}
Assume $p=(s, T) \in \Pp$ and $D\subseteq\Pp$ is dense open. Then   there are $q=(s, S)\leq^*p$
and $\ell\geq\lh(s)$ such that
\[
 (\forall t\in\Lev_\ell(S))\quad (t, S\restr t)\in D.
\]
\end{lemma}
\begin{proof}
The proof is as usual, we present it for completeness.
For $t\in T$ extending $s$, let us say $t$ is good if there exist $S_t$ and $r_t<\omega$
   such that $q=(t, S_t)\leq^*p[t]$ and
$q[u] \in D$ for every
$u\in\Lev_{\lh(t)+r_t}(S_t)$. It is easily seen that if
 $t$ is not good, then
\[
 B_t=\{\alpha\in\Succ_T(t):t\concat\langle\alpha\rangle
                         \text{ is good}\}\in I(t).
\]
To see this, suppose by the way of contradiction that $B_t$ is positive. For each $\alpha \in B_t$
let $S_\alpha$ and $r_\alpha$ witness $t\concat\langle\alpha\rangle$ is good. Since $I(t)^+$ is countably complete, we can find
a positive subset $B'_t\subseteq B_t$ and $r< \omega$ such that $r_\alpha=r$ for all $\alpha \in B'_t$.
Then $S=\bigcup_{\alpha \in B'_t}S_\alpha$ and $r+1$ witness $t$ is good which is a contradiction.

To prove the lemma, it suffices to show that $s$ is good. If not, using the above observation, we can recursively thin out the tree $T$
to get a subtree $S$ of $T$ such that every $t \in S$ is not good. Let $q=(s, S)$, and choose $r \leq q$ such that $r \in D$.  Set $t=s^r \in T.$
Then $T^r$ and $r_t=0$ witness that $t$ is good, a contradiction.
\end{proof}
By the strong Prikry property, we have the following.
\begin{corollary}[Prikry property]\label{cor:pure}
For every forcing sentence $\varphi$ and $p \in \Pp$, there exists
$q\leq^*p$ deciding $\varphi$.
\end{corollary}
The following is also immediate.
\begin{lemma}
Assume $p, q \in \Pp$ and $s^p=s^q$. Then $p$ and $q$ are compatible. In particular $\Pp$ satisfies the $\kappa^+$-c.c.
\end{lemma}

\begin{lemma}\label{lem:tail}
Assume $p=(s, T) \in \Pp$ and $\lh(s)=m$. Then
$
 (\{q \in \Pp: q\leq^*p\},\leq^*)
$
is $\lambda_m$-strategically closed.
\end{lemma}
\begin{proof}
This follows easily from the fact that if $t \in T$ has
height $n\geq m$, then $(I(t)^+, \subseteq)$ is $\lambda_m$-strategically closed.
\end{proof}
Using Corollary~\ref{cor:pure} an
Lemma~\ref{lem:tail}, we get the following.
\begin{corollary}\label{cor:short}
Assume $p=(s, T) \in \Pp$, $\lh(s)=m$, $|c|<\lambda_m$ and
 $p\forces\dot a\subseteq\check c$. Then there exist
$q\leq^*p$ and $a\subseteq c$ such that
$q\forces\dot a=\check a$.
\end{corollary}
By Corollary~\ref{cor:short}, we have the following.
\begin{theorem}\label{thm:bounded}
forcing with $\Pp$ adds no new bounded subsets of $\kappa$.
\end{theorem}
Putting the above results together we get the following.

\begin{lemma}\label{lem:regularpreservation}
Forcing with $\Pp$ preserves all cardinals and cofinalities.
\end{lemma}

We now study some further basic properties of the forcing notion $\Pp.$ We start with the following easy lemma.
\begin{lemma}\label{lem:fusion}
Assume $\langle p_j:j<\omega\rangle$ is a decreasing sequence of conditions with a common stem $s$.
Assume  there are integers $\langle h_j: j< \omega \rangle$ such that
$ T^{p_{j+1}}\cap{}^{\leq h_j}\Ord
       =T^{p_j}\cap{}^{\leq h_j}\Ord$
for $j<\omega$, and $\lim_{j<\omega}h_j=\infty$. Set
 $T=\bigcap_{j<\omega}T^{p_j}$. Then $q=(s, T)$ is a condition.
\end{lemma}

The next theorem shows that any set of ordinals of size $\kappa$ in the extension is coversed by a ground model set of size $\kappa$.
\begin{theorem}\label{thm:cover}
Assume  $p\forces\dot f:\check\kappa\to\Ord$.
Then there are $q\leq^*p$ and a set $B$ of ordinals such that $|B| \leq \kappa$ and
$q\forces\ran(\dot f)\subseteq\check B.$
\end{theorem}
\begin{proof}
Let $\langle m_j: j<\omega   \rangle$ be an increasing integers  of integers such that
 $m_j\geq\lh(s^p)$ and
$\kappa_j<\lambda_{m_j}$. We define a decreasing sequence $\langle p_j:j<\omega\rangle$ of direct extensions of $p$ as in Lemma~\ref{lem:fusion}.
Set $p_0=p$.
Assume now $j<\omega$ and  $p_j$ is defined. Fix
$t\in\Lev_{m_j}(T^{p_j})$.
For $\alpha<\kappa_j$, let $D^j_{t,\alpha}$ consist of the pure extensions
$q\leq^*p_j[t]$ such that for some $\ell\geq\lh(t)$, $q[t]$ decides
$\dot f(\alpha)$ for all  $t\in\Lev_\ell(T^q)$.
By Lemma~\ref{lem:strong}, $D^j_{t,\alpha}$ is $\leq^*$-dense open. As $\lambda_{m_j} > \kappa_j$, it follows from
 Lemmas~\ref{lem:tail} that $\bigcap_{\alpha<\kappa_j}D^j_{t,\alpha}$ is dense open, so let $q_{j, t}$ be in this intersection.
 Then for each $\alpha<\kappa_j$, we can find $\ell(t, j, \alpha)\geq m_j$ and ordinals $\langle  \xi_{t, j,  s, \alpha}: s\in\Lev_{\ell(t, j, \alpha)}(T^{q_{j, t}})    \rangle$ such that $q_{j, t}[s] \forces \dot f(\alpha)=\check \xi_{t, j, s, \alpha}$.
 Let $p_{j+1}=(s, T^{p_{j+1}})$, where
  $T^{p_{j+1}}=\bigcup_{t\in\Lev_{m_j}(T^{p_j})}T^{q_{j, t}}$. The hypotheses of Lemma~\ref{lem:fusion} are satisfied, hence $q=(s, T)$
  is a condition where $T=\bigcap_{j<\omega}T^{p_j}$. Let also
  \[
  B=\{\xi_{t, j, s, \alpha}: j<\omega,~t\in\Lev_{m_j}(T^{p_j}), \alpha<\kappa_j \text{ and } s\in\Lev_{\ell(t, j, \alpha)}(T^{q_t})    \}.
  \]
The set $B$ has clearly cardinality $\leq \kappa$, and by our construction, it is easy to see that $q\forces\ran(\dot f)\subseteq\check B.$
\end{proof}

Assume $G$ is $\Pp$-generic over $V$ and set $C=\bigcup\{s^p:p\in G\}$. Then, as in the usual Prikry forcing, we have the following.
\begin{lemma}\label{lem:branch}
\begin{enumerate}
\item $C\in\prod_{n<\omega}\kappa_n$ and $C(n)\in(\kappa_{n-1},\kappa_n)$, for all $n>0$,
\item $G=\{p \in \Pp:C \in[T^p]\}$, in particular $V[G]=V[C].$
\end{enumerate}
\end{lemma}

\subsection{Some canonization results}
We now prove a few results about our ideals, that will be used later.
\begin{lemma}\label{lem:onecoordinate}
Assume $I$ is a normal ideal on a regular uncountable cardinal $\rho$,
 $A\in I^+$, and  $f:A\to Y$.
Then there is $B\subseteq A$ such that $B\in I^+$ and $f\restr B$ is constant
or injective.
\end{lemma}
\begin{proof}
If $f^{-1}\{y\}\in I^+$ for some $y\in Y$, then we are done.  So assume that $f^{-1}\{y\}\in I$ for each $y\in Y$.
Let
\[
 C=\{\alpha\in A:(\exists\beta\in A\cap\alpha)\ f(\beta)=f(\alpha)\}.
\]
We claim that $C \in I$. Assume not. Then
$C\in I^+$. Define the function $g$ on $C$ by
\[
 g(\alpha)=\min\{\beta\in A\cap\alpha:f(\beta)=f(\alpha)\}.
\]
Then $g$ is regressive on $C$, hence by Lemma~\ref{lem:normalfacts}(iii), there are $C'\subseteq C$ in $I^+$
and $\beta$ such that $g\restr C'$ is constantly equal to $\beta$.
But then $C'\subseteq f^{-1}\{f(\beta)\}$, a contradiction to $f^{-1}\{f(\beta) \in I$.
Hence $C\in I$, and $f$ is injective on $A\setminus C\in I^+$.
\end{proof}

\begin{lemma}\label{lem:separation}
Assume $\tau<\rho$,  $I_i$ is a normal ideal on $\rho$ for $i<\tau$,
and  $A_i\in I_i^+$ are pairwise disjoint.
Assume that  $f_i:A_i\to Y$ is injective for every $i<\tau$. Then there are sets
$B_i\subseteq A_i$ such that $B_i\in I_i^+$ and the sets $f_i[B_i]$
are pairwise disjoint.
\end{lemma}
\begin{proof}
For  $i\neq j<\tau$, set
\[
 E_{i,j}=\{\alpha\in A_i:(\exists\beta\in A_j\cap\alpha)\
                                      f_i(\alpha)=f_j(\beta)\}.
\]
Since $f_j$ is injective, for each $\alpha\in E_{i,j}$, there is a unique witness $\beta$ as above, call it $g_{i, j}(\alpha)$.
 The function $g_{i, j}$ is
regressive and, since $f_i$ is injective, it is injective as well. By Lemma~\ref{lem:normalfacts},
 $E_{i,j}\in I_i$, and hence $\bigcup_{j\ne i}E_{i,j} \in I_i$.
Thus
$
 B_i=A_i\setminus\bigcup_{j\ne i}E_{i,j} \in I_i^+.
$
It is easily seen that the sets $B_i$  for $i<\tau$ are as required.
\end{proof}
As corollary, we have the following finite-level canonization theorem.
\begin{theorem}\label{thm:canon}
Assume $p \in \Pp$, $m=\lh(s^p)$, $m\leq\ell<\omega$, and
$f:\Lev_\ell(T^p)\to Y$.
Then there are $q\leq^*p$, $k\in[m,\ell]$, and an injective function
$h:\Lev_k(T^q)\to Y$ such that for all $u\in\Lev_\ell(T^q)$, $f(u)=h(u\restr k).$
\end{theorem}
\begin{proof}
We prove the theorem by induction on
 $\ell \geq m$. For simplicity set $p=(s, T)$. If $\ell=m$, this is immediate as $\Lev_m(T)=\{s\}$ is singleton.
Suppose now that $\ell>m$. For each
$t\in\Lev_{\ell-1}(T)$ define $f_t: \Succ_T(t)  \rightarrow Y$ by
$
 f_t(\alpha)= f(t\concat\langle\alpha\rangle).
$
By Lemma~\ref{lem:onecoordinate}, there is a positive set $A_t \subseteq \Succ_T(t)$ on which $f_t$ is
constant or injective. Now define a coloring $c: \Lev_{\ell-1}(T) \rightarrow \{0, 1\}$
so that $c(t)=0$ iff $f_t \restr A_t $ is constant (and hence $c(t)=1$ iff $f_t \restr A_t $ is injective). By
Lemma~\ref{lem:homogeneous}, there is $q \leq^* p$ such that $c \restr \Lev_{\ell-1}(T^q)$ is constant. We now consider two case.

First assume that  $c \restr \Lev_{\ell-1}(T^q)$ is equal to $0$.
Define $g: \Lev_{\ell-1}(T^q) \to Y$ by $g(t)=f(t\concat\langle\alpha\rangle)$, where
 $\alpha \in \Succ_{T^q}(t)$. Note that the function $g$ is well-defined and its definition is independent of the choice of $\alpha$.
By
the induction hypothesis applied to $g$ at height $\ell-1$, we get some $r \leq^* q$,
 $k\leq\ell-1$, and $h$.
Then, for  $u\in\Lev_\ell(T^r)$, we have
$
 f(u)=g(u\restr(\ell-1))=h(u\restr k).
$

Next suppose that  $c \restr \Lev_{\ell-1}(T^q)$ is equal to $1$.
Now the sets $\langle A_t: t\in   \Lev_{\ell-1}(T^q)  \rangle$ are pairwise disjoint and
each $f_t \restr A_t $ is injective. Thus by
 Lemma~\ref{lem:separation}, we can find  positive sets $B_t \subseteq A_t$
 such that the sets $f_t[B_t]$ are pairwise disjoint.
The resulting pure refinement $r$ makes $f$ injective on
$\Lev_\ell(T^r)$. Take $k=\ell$ and $h=f\restr\Lev_\ell(T^r)$.
\end{proof}

\subsection{Minimality for subsets of $\kappa$}
The next lemma follows easily from Corollary~\ref{cor:short} and Lemma~\ref{lem:fusion}.
\begin{lemma}\label{lem:reading}
Assume $p \in \Pp$,  $p\forces\dot A\subseteq\check\kappa$ and
$m=\lh(s^p)\geq1$.
Then there are $q\leq^*p$ and a family
$
 \langle a_t:t\in T^q,\ \lh(t)\geq m\rangle
$
such that, whenever $\lh(t)=k\geq m$,
$
 a_t\subseteq\kappa_{k-1}$ and
 $q[t]\forces\dot A\cap\check\kappa_{k-1}=\check a_t.$
Furthermore, if $u\subseteq t$ and $\lh(u)\geq m$, then
$a_t\cap\kappa_{\lh(u)-1}=a_u$.
\end{lemma}
\begin{lemma}\label{lem:localcoding}
Assume $p$ and the sequence $\langle a_t\rangle$ are as in
Lemma~\ref{lem:reading}, and
$p\forces\dot A\notin\check{\pow(\kappa)^V}.$
Then there are $q\leq^*p$, an integer $k>\lh(s^p)$, an ordinal
$\beta<\kappa$, and an injective family
$\langle d_t:t\in\Lev_k(T^q)\rangle$ of subsets of $\beta$ such that for every $t\in\Lev_k(T^q)$,
$
q[t]\forces\dot A\cap\check\beta=\check d_t.
$
\end{lemma}
\begin{proof}
Let $m=\lh(s^p)$. By Lemma~\ref{lem:tail}, fix
 a winning strategy for player II  in the  $\omega$-length game
 on pure extensions of $p$.
At step $j<\omega$, let $p_j$ be the strategy played by II and set
$\ell_j=m+j+1$. Define the map $f_j: \Lev_{\ell_j}(T^{p_j}) \to  \pow(\kappa)$ by
$f_j(u)=a_u$.
By Theorem~\ref{thm:canon}, there are
 $r_j\leq^*p_j$, a $k_j\in[m,\ell_j]$, and an injective
function $h_j$ as there.

We show that there is $j<\omega$ such that $k_j>m$. Suppose not.
Then we run the above game, where after the choice of $p_j$, player I plays
$r_j$. Note that for each $j$,
$r_j\forces \dot A\cap\check\kappa_{\ell_j-1}=\check z_j$, for some $z_j$ and the sequence $\langle z_j: j<\omega \rangle$
is coherent. Let $q_\omega$
be the condition played by II at step $\omega$. It then extends each $r_j$,
and hence
$q_\omega \forces \dot A=\check{\bigcup_{j<\omega}z_j}$, which is a contradiction.

Thus let $j<\omega$ be such that $k_j>m$. Then $q=r_j,$  $k=k_j$, and
$\beta=\kappa_{\ell_j-1}$ are as required.
\end{proof}

\begin{lemma}\label{lem:codingsystem}
Assume $p \in \Pp,$
$p\forces$``$\dot A\subseteq\check\kappa$ and
$\dot A\notin\check{\pow(\kappa)^V}$''.
Then there are $q\leq p$, fronts $\langle F_j:j<\omega\rangle$ of $q$,
and, for each $t\in F_j$, an ordinal $\beta_t<\kappa$ and a sequence
$
 \langle d_{t,u}:u\in F_{j+1},\ t\subsetneq u\rangle
$
such that the following conditions are hold.
\begin{enumerate}
\item The fronts are $\prec$-increasing, $F_0=\{s^q\}$,
and every member of $F_j$ has height at least $\lh(s^q)+j$.
\item If $t\in F_j$, then the sets $d_{t,u}\subseteq\beta_t$ are
pairwise distinct.
\item If $t\in F_j$, $u\in F_{j+1}$, and $t\subsetneq u$, then $q[u]\forces\dot A\cap\check\beta_t=\check d_{t,u}$.
\end{enumerate}
\end{lemma}
\begin{proof}
By induction on $j<\omega$ we define a $\leq^*$-decreasing sequence $\langle p_j: j<\omega   \rangle$ of conditions and a $\prec$-increasing sequence
of fronts as follows.
Let $p_0 \leq^* p$ satisfies the conclusion of
Lemma~\ref{lem:reading}, and set  $F_0=\{s^p\}$.
Now assume that $j<\omega$ and
 $p_j$ and $F_j$ are constructed.

 For each $t \in F_j$, by Lemma~\ref{lem:localcoding}, we can find some $q_t \leq^* p_j[t]$,
 $k_t>\lh(t)$, and ordinal $\beta_t$, and an injective family
$\langle d_{t,u}:u\in\Lev_{k_t}(T^{q_t})\rangle$ as there. Let $p_{j+1}=(s^p, T^{p_{j+1}})$, where $T^{p_{j+1}}=\bigcup_{t\in F_j}T^{q_t}$
 and $F_{j+1}=\bigcup_{t\in F_j}\Lev_{k_t}(T^{q_t}).$
 By Lemma~\ref{lem:gluing}, $p_{j+1}$ is a condition and $F_{j+1}$ is a front, and clearly $F_j \prec F_{j+1}$.
By Lemma~\ref{lem:fusion}, $q=(s, T^q)$ is a condition, where $T^q=\bigcap_jT^{p_j}$. It is clear that
 $q$, fronts $\langle F_j:j<\omega\rangle$,
 ordinals $\beta_t<\kappa$ and
$
 \langle d_{t,u}:u\in F_{j+1},\ t\subsetneq u\rangle,
$
where $t\in F_j$ are as required.
\end{proof}
Finally we can prove our minimality result for subsets of $\kappa.$
\begin{theorem}\label{thm:subsetminimal}
Assume $G$ is $\Pp$-generic over $V$ and
$A\in\pow(\kappa)^{V[G]}\setminus V$, then $V[A]=V[G]$.
\end{theorem}
\begin{proof}
Let  $\dot A$ be a $\Pp$-name for $A$ and let  $p \in G$ forces
$\dot A\subseteq\check\kappa$ and
$\dot A\notin\check{\pow(\kappa)^V}$.
By Lemma~\ref{lem:codingsystem}
we can find  $q \in G$, fronts $\langle F_j:j<\omega\rangle$ of $q$, and for
$t\in F_j$, an ordinal $\beta_t<\kappa$ and the sequence
$
 \langle d_{t,u}:u\in F_{j+1},\ t\subsetneq u\rangle,
$
as there.
Let also $C=\bigcup\{s^p: p \in G    \}$ be the generic branch.

Working in $V[A]$, we define by recursion a sequence  $\langle t_j:j<\omega\rangle$ with $t_j \in F_j$ and  $t_j\subset C$   as follows.
Let $t_0=s^q$.
Now assume that $j<\omega$ and $t_j\in F_j$ is chosen so that $t_j\subset C$. Let $t_{j+1}$ be the unique $u\in F_{j+1}$
which properly extends $t_j$ and satisfies
$
 A\cap\beta_{t_j}=d_{t_j,u}.
$
Since $\lh(t_j)\geq\lh(s^q)+j$,
$C=\bigcup_{j<\omega}t_j\in V[A]$. Hence by
Lemma~\ref{lem:branch},
$V[G]=V[C]\subseteq V[A]\subseteq V[G]$.
\end{proof}
The following corollary immediate using Theorem~\ref{thm:cover}.
\begin{corollary}\label{cor:smallminimal}
Assume $X\in V[G]\setminus V$ is a set of ordinals and
$V[G]\models |X|\leq\kappa$, then $V[X]=V[G]$.
\end{corollary}

\subsection{Minimality for arbitrary sets of ordinals}
We now show that the above minimality argument can be extended to all sets of ordinals. For this we prove a general result which reduces minimality for arbitrary sets of ordinals to subsets of $\kappa$.
\begin{lemma}\label{lem:subsetreduction}
Assume $\kappa$ is an infinite cardinal and $\Pp$ is a forcing
notion satisfying the following conditions:
\begin{enumerate}
\item $\Pp$ is $\kappa^+$-c.c.
\item There are a set $S$ and a map $\sigma:\Pp\to S$ such that
$|S|\leq\kappa$ and, if $\sigma(p)=\sigma(q)$, then
$p$ and $q$ are compatible.
\item For every $\Pp$-generic filter $G$ over $V$ and every set
$A\in\pow(\kappa)^{V[G]}\setminus V$, one has $V[A]=V[G]$.
\end{enumerate}
Then, for every $\Pp$-generic filter $G$ over $V$ and every set of
ordinals $X\in V[G]\setminus V$,  $V[X]=V[G]$.
\end{lemma}

\begin{proof}
Let $G$ be a  $\Pp$-generic filter  over $V$ and let $X=\dot X[G] \in V[G]\setminus V$ be  a set of ordinals. Let $\theta$ be the least ordinal such that
$
 a=X\cap\theta\notin V.
$
By minimality, $\theta$ is a limit ordinal,  $a$ is unbounded in $\theta$
and $a\cap\beta=X\cap\beta\in V$    for every $\beta<\theta$.
Let $\lambda=\cf^V(\theta)$ and let $\langle\theta_i:i<\lambda\rangle \in V$  be an increasing
and cofinal sequence  in $\theta$.
For each $i<\lambda$, set
$
 a_i=a\cap\theta_i=X\cap\theta_i\in V.
$
For each $i<\lambda$, let $p_i\in G$ be such that
$
 p_i\forces\dot X\cap\check{\theta_i}=\check{a_i}.
$
We claim that $\lambda\leq\kappa$. Suppose not. Then
as $|S|\leq\kappa<\lambda$, by clause (2), there is
$s\in S$ such that
$
 I=\{i<\lambda:\sigma(p_i)=s\}
$
is unbounded in $\lambda$.
In $V$, define the relation $h_s$ as follows:
\[
 \begin{split}
 h_s=\{(i,b)\in\lambda\times\pow(\theta)^V:
       &(\exists q\in\Pp)\bigl(\sigma(q)=s\text{ and }
      q\forces
                  \dot X\cap\check{\theta_i}=\check b\bigr)\}.
 \end{split}
\]
Clearly $I\subseteq\dom(h_s)$.
By compatibility assumption in (2),  $h_s$ is a  partial function.
Let $i\in\dom(h_s)$, $i<j \in I$ and let $q$ witness the value $h_s(i)$.
The conditions $q$ and $p_j$ have the same $\sigma$-value, hence they are compatible and
their common extension forces
$
 \check{h_s(i)}
       =\dot X\cap\check{\theta_i}
       =\check{a_j}\cap\check{\theta_i}.
$
Thus
$
 h_s(i)=a_j\cap\theta_i=a_i.
$
Since $I$ is unbounded in $\lambda$ and
$I\subseteq\dom(h_s)$, it follows that
\[
 a=\bigcup_{i\in I}a_i
   =\bigcup_{i\in\dom(h_s)}h_s(i)
   =\bigcup\ran(h_s)\in V,
\]
which contradicts the choice of $a$.  Thus $\lambda\leq\kappa$.

If $\theta\leq\kappa$, then the conclusion follows from clause (3),
thus suppose that $\theta>\kappa$.
Set $D=\{(i,b):i<\lambda,\ b\in D_i\}$, where
for $i<\lambda$,
\[
 D_i=\{b\in\pow(\theta)^V:(\exists q\in\Pp)\
                       q\forces\dot X\cap\check{\theta_i}=\check b\}.
\]
Then $D \in V$.
For each $i$, by clause (1), $|D_i| \leq \kappa$ and hence $|D| \leq \kappa$.
Working in $V[G]$ set $Y=\{(i,a_i):i<\lambda\}$. By the choice of conditions $p_i$, we have $Y\subseteq D$.
It then follows from clause (3) that
$
V[a]=V[Y]=V[G],
$
and hence $V[X]=V[G]$ as required.
\end{proof}
Note that our forcing notion $\Pp$ satisfies items (1)--(3) of the above lemma, and hence its conclusion can be applied to it.
Putting the above results together, we get the following.

\begin{theorem}\label{thm:fullminimal}
Assume $G$ is a $\Pp$-generic filter over $V$. Then
for every set of ordinals
$X\in V[G]\setminus V$,  $V[X]=V[G]$.
\end{theorem}

\subsection{Proof of the main theorem}
We are finally ready to prove our first main theorem.

\begin{proof}[Proof of Theorem~\ref{thm:main}]
Assume the sequence $\langle \kappa_n: n< \omega  \rangle$ of measurable cardinals is as in
 Theorem~\ref{thm:main}. Let $V_1$ be the generic extension obtained by Theorem~\ref{thm:preparation}. In $V_1$, force with the forcing notion
 of Section 2 and let $V_2$ be the resulting generic extension. Then $V_1$ and $V_2$ are as required.
\end{proof}

\section{Forcing a minimal subset of \texorpdfstring{$\aleph_\omega$}{aleph omega} using one measurable cardinal}\label{sec:interleaved}
In this section we introduce another method for adding a minimal fresh subset to $\aleph_\omega$, using just one measurable cardinal.
This will give us an equi-consistency result.

We start with a model in which $\kappa$ is a measurable cardinal and there exists a sequence $\langle U_\alpha:\alpha<\kappa\rangle$
 of distinct normal measures on $\kappa$. Let
$\langle X_\alpha:\alpha<\kappa\rangle$ be pairwise disjoint subsets of $\kappa$ such that
$X_\alpha\in U_\alpha$, for all $\alpha < \kappa$.  By shrinking the measure one sets, we may assume that
each set $X_\alpha$ consists of inaccessible cardinals above $\alpha$.

Set  $U=U_0$, and let
$ j:V\longrightarrow M=\Ult(V,U)$ be the corresponding elementary embedding with critical point $\kappa$ and set $\kappa_1=j(\kappa)$.
For $0<\alpha<\kappa$ let
$  i_\alpha:M\longrightarrow N_\alpha  =\Ult(M,j(U_\alpha))$
       be the two step iteration and set
$ k_\alpha=i_\alpha\circ j$  and
 $ \kappa_{2,\alpha}=i_\alpha(\kappa_1).$
It is easy to see there are $K_1$ and $K_{2, \alpha,},  0<\alpha<\kappa$, in $V$, such that
\begin{itemize}
\item $K_1$ is $\Col(\kappa^+,<\kappa_1)^M$-generic over $M$,

\item $K_{2, \alpha}$ is $\Col(\kappa^+,[\kappa_1,\kappa_{2,\alpha}))^{N_\alpha}$-generic over $N_\alpha[K_1].$
\end{itemize}

\subsection{The forcing notion and its basic properties}
In this subsection, we introduce our main forcing notion and prove its basic properties.

\begin{definition}\label{ic:def:P}
A condition in $\Pp$ is of the form
\[
 p=\bigl(\langle(\nu^p_i,\xi^p_i,f^p_i):i<m^p\rangle,
 A^p,\langle B^p_\alpha:0<\alpha<\kappa\rangle,F^p_1,F^p_2\bigr),
\]
where:
\begin{enumerate}
\item $0<m^p=m<\omega$ and $\nu^p_0<\xi^p_0<\nu^p_1<\xi^p_1<\cdots<\nu^p_{m-1}<\xi^p_{m-1}<\kappa$;
\item $\xi^p_i\in X_{\xi^p_{i-1}}$ for $0<i<m$;
\item $f^p_i\in\Col((\nu^p_i)^+,<\nu^p_{i+1})$ for $i<m-1$, and
$f^p_{m-1}\in\Col((\nu^p_{m-1})^+,<\kappa)$;
\item $A^p \subseteq X_0$ is in $U$ and
$\min(A^p)>\max\{\xi^p_{m-1},\sup(\supp(f^p_{m-1}))\}$;
\item $B^p_\alpha \subseteq X_\alpha$ is in $U_\alpha$ for $0<\alpha<\kappa$;
\item $F^p_1(\nu)\in\Col(\nu^+,<\kappa)$ for $\nu\in A^p$;
\item $ j(F^p_1)(\kappa)\in K_1;$
\item $F^p_2$  has domain $T_p=\{(\nu,\xi):\nu\in A^p,\ \xi\in B^p_\alpha,  \ 0<\alpha<\nu<\xi\}$
and for every $(\nu,\xi) \in T_p,$ $F^p_2(\nu,\xi)\in\Col(\nu^+,[\xi,\kappa))$;
\item $ k_\alpha(F^p_2)(\kappa,\kappa_1)\in K_{2,\alpha}$ for
$0<\alpha<\kappa.$
\end{enumerate}

For conditions $p, q \in \Pp$, we say $q \leq p$ if
\begin{itemize}
\item  $m^p \leq m^q$ and for all $i< m^p,$ $(\nu^p_i,\xi^p_i)=(\nu^q_i,\xi^q_i)$ and $f^q_i \leq f^p_i$,
\item For $m^p \leq i < m^q, \nu^q_i \in A^p$,
$\xi^q_i \in \bigcup_{\alpha>0}B^p_\alpha$ and $\sup(\supp(F^p_1(\nu^q_i))) < \xi^q_i$,

\item For $m^p \leq i < m^q, f^q_i \leq F^p_1(\nu^q_i) \cup F^p_2(\nu^q_i, \xi^q_i)$,
\item $A^q \subseteq A^p,$ and $B^q_\alpha \subseteq B^p_\alpha$ for $0<\alpha < \kappa,$
\item  For $\nu \in A^q, F_1^q(\nu) \leq F^p_1(\nu)$,

\item If $(\nu, \xi) \in T_q$,
 then $F^q_2(\nu, \xi) \leq F^p_2(\nu, \xi)$.
\end{itemize}
We also set $q \leq^* p$ if $q \leq p$ and $m^q=m^p$.
\end{definition}
\begin{notation}
Assume $p \in \Pp$ is of length $m=m^p$. Then  we set
\begin{itemize}
\item  $\stem(p)=\langle(\nu^p_i,\xi^p_i,f^p_i):i<m\rangle$, and call it the stem of $p$,

\item $\tail(p)=( A^p,\langle B^p_\alpha:0<\alpha<\kappa\rangle,F^p_1,F^p_2)$ and call it the tail of $p$.

\item $\vec\nu_p=\langle  \nu^p_0, \nu^p_1, \cdots, \nu^p_{m-1} \rangle$,

\item $\vec\xi_p=\langle  \xi^p_0, \xi^p_1, \cdots, \xi^p_{m-1} \rangle$,

\item $\vec f_p=\langle  f^p_0, f^p_1, \cdots, f^p_{m-1} \rangle$.
\end{itemize}
\end{notation}
We have the following easy lemma.
\begin{lemma}\label{ic:normal}
Assume $p \in \Pp$. Then there exists $q \leq^* p$ such that:
\begin{enumerate}
\item If $\nu\in A^q$, $\xi\in\bigcup_{\alpha>0}B^q_\alpha$,
and $\nu<\xi$, then $\sup(\supp(F^q_1(\nu)))<\xi$.
\item If $\nu\in A^q$, $\xi,\zeta\in\bigcup_{\alpha>0}B^q_\alpha$,
and $\nu<\xi<\zeta$, then $\sup(\supp(F^q_2(\nu,\xi)))<\zeta$.
\item If $\nu<\lambda$ are in $A^q$, then
$\sup(\supp(F^q_1(\nu)))<\lambda$.
\item If $\nu,\lambda\in A^q$,
$\xi\in\bigcup_{\alpha>0}B^q_\alpha$, and $\nu<\xi<\lambda$,
then $\sup(\supp(F^q_2(\nu,\xi)))<\lambda$.
\item If $\lambda\in A^q$ and $0<\alpha<\lambda$, then
$\sup(B^q_\alpha\cap\lambda)=\lambda$.
\end{enumerate}
\end{lemma}

\begin{definition}
\begin{enumerate}
\item We say a condition $q$ is \emph{normal} when it satisfies the conclusion of  Lemma~\ref{ic:normal}. By the lemma, the set of
 normal conditions is $\leq^*$-dense in $\Pp$.


 \item Assume $p \in \Pp$ and $t=\langle (\nu_i, \xi_i): m^p \leq i< n    \rangle$ of length $n-m^p$ is such that $\xi^p_{m^p-1} < \nu_{m^p} < \xi_{m^p} < \cdots < \nu_{n-1} < \xi_{n-1}$, $\nu_i \in A^p$, $\xi_i \in \bigcup_{\alpha>0}B^p_\alpha$ and $\sup(\supp(F^p_1(\nu_i)))< \xi_i$. Then we define the condition $p[t]$ of length $n=\lh(p)+\lh(t)$ as
\begin{itemize}
\item $\stem(p[t])= \bigl(\langle(\nu^p_i,\xi^p_i,f^p_i):i<m^p\rangle^{\frown} \langle (\nu_i, \xi_i, F^p_1(\nu_i) \cup F^p_2(\nu_i, \xi_i)): m^p \leq i<n  \rangle,$

\item $\tail(p[t])=( A^p \cap (\xi_{n-1}, \kappa),\langle B^p_\alpha \cap (\xi_{n-1}, \kappa) :0<\alpha<\kappa\rangle,F^{p[t]}_1,F^{p[t]}_2\bigr),$
where $F^{p[t]}_1,F^{p[t]}_2$ are the restriction of $F^{p}_1,F^{p}_2$ on their new domain respectively.
\end{itemize}
In this case we say $t$ is legal for $p$.
 \end{enumerate}
\end{definition}

Assume $G$ is $\Pp$-generic over $V$ and set
\begin{itemize}
\item  $\vec\kappa=\langle\kappa_n:n<\omega\rangle,$ where $\kappa_n=\nu_n^p$, for some (and hence all) $p \in G$ with $m^p > n$,
\item $\vec\eta=\langle\eta_n:n<\omega\rangle,$ where $\eta_n=\xi_n^p$, for some  $p \in G$ with $m^p > n$,
\item $\vec H=\langle H_n:n<\omega\rangle,$ where $H_n=\{f^p_{n}: p \in G, m^p > n    \}$. Then each $H_n$ is $\Col(\kappa_n^+,<\kappa_{n+1})$-generic over $V$.
\end{itemize}
As in the usual Prikry forcing, we have the following.
\begin{lemma}
\begin{enumerate}
\item The sequence $\langle\kappa_n:n<\omega\rangle$ is cofinal in $\kappa,$
\item For every $A \in U$, there exists $n<\omega$ such that $\kappa_i\in A$, for all $i \geq n.$
\end{enumerate}
\end{lemma}
We also have the Prikry property.
\begin{lemma}
\label{prikry property}

\begin{enumerate}
\item Assume $p \in \Pp$ and $D \subseteq \Pp$ is dense open. Then there are $q \leq^* p$ and $\ell$ such that if $t$ is legal for $q$ and $\lh(t)=\ell$, then $q[t] \in D$.

\item The forcing notion $(\Pp, \leq, \leq^*)$ satisfies the Prikry property.
\end{enumerate}
\end{lemma}

\begin{definition}
Assume $p \in \Pp$ and $n < m^p$. Then
\begin{itemize}
\item $p \res n= \langle(\nu^p_i,\xi^p_i,f^p_i):i<n\rangle,$
\item $p\setminus n=\bigl(\langle(\nu^p_i,\xi^p_i,f^p_i):n \leq i< m^p\rangle,  A^p,\langle B^p_\alpha:0<\alpha<\kappa\rangle,F^p_1,F^p_2\bigr),$
\end{itemize}
\end{definition}

We have the following factorization
property.
\begin{lemma}\label{lem:factor}
Assume $p \in \Pp$ and $n < m^p$.
Then
\[
\Pp/p \simeq \prod_{i< n} \Col((\nu_i^p)^+, < \nu^p_{i+1}) \times \Pp/(p\setminus n).
\]
Furthermore,
\begin{enumerate}
\item $\prod_{i< n} \Col((\nu_i^p)^+, < \nu^p_{i+1})$ is $\nu^p_{n}$-c.c.,

\item $(\Pp/(p\setminus n), \leq^*)$ is  $(\nu^p_{n})^+$-closed.

\item For every $\beta \leq \kappa_{n}$,  $\mathcal{P}(\beta)^{V[G]}=\mathcal{P}(\beta)^{V[\langle H_i: i< n  \rangle]}$.
\end{enumerate}
\end{lemma}
The following lemma is also immediate.
\begin{lemma}
\label{lem:chain} The forcing $(\Pp, \leq)$ satisfies the $\kappa^+$-c.c. Indeed if $p, q \in \Pp$ have the same stem, then they are compatible.
\end{lemma}

Set
$V_1=V[\vec\kappa,\vec H],$ and $V_2=V_1[\vec\eta]\subseteq V[G].$
We will show that the pair $(V_1,V_2)$ has the required properties.

We first determine the forcing notion that leads to the generic extension $V_1$.
\begin{definition}\label{ic:def:B}
A condition in $\Q$ if of the form
\[
 b=(\langle(\nu^b_i,f^b_i):i<m^b\rangle,A^b,F^b),
\]
where
\begin{enumerate}
\item $m^b<\omega$ and $\nu^b_0<\nu^b_1<\cdots<\nu^b_{m-1}<\kappa$;
\item $f^b_i\in\Col((\nu^b_i)^+,<\nu^b_{i+1})$ for $i<m-1$, and
$f^b_{m-1}\in\Col((\nu^b_{m-1})^+,<\kappa)$;
\item $A^b \subseteq X_0$ is in $U$ and
$\min(A^b)>\sup(\supp(f^b_{m-1}))$;
\item $F^b_1(\nu)\in\Col(\nu^+,<\kappa)$ for $\nu\in A^b$.
\end{enumerate}
We set $\stem(b)=\langle(\nu^b_i,,f^b_i):i<m^b\rangle$ and call it the stem of $b$,  similarly, set
$\tail(b)=( A^b,F^b_1)$ and call it the tail of $b$.

For conditions $b, c \in \Q$, we say $c \leq b$ if
\begin{itemize}
\item  $m^b \leq m^c$ and for all $i< m^b,$ $\nu^b_i=\nu^c_i$ and $f^c_i \leq f^b_i$,
\item For $m_b \leq i < m^c, \nu^b_i \in A^b$,

\item For $m_b \leq i < m^c, f^c_i \leq F^c_1(\nu^c_i)$,
\item $A^c \subseteq A^b,$
\item  For $\nu \in A^c, F_1^c(\nu) \leq F^b_1(\nu)$.
\end{itemize}
We also set $c \leq^* b$ if $c \leq b$ and $m^c=m^b$.
\end{definition}
Define $\pi: \Pp \to \Q$ by
\[
 \pi(p)=(\langle(\nu_i^p,f_i^p):i<m^p\rangle,A^p,F_1^p).
\]
We have the following easy lemma

\begin{lemma}\label{ic:projection}
The map $\pi: \Pp \to \Q$ is a projection, i.e., it is order preserving, and if $b\leq\pi(p)$, there is $q\leq p$ in $\Pp$ such that
$\pi(q)\leq b$.
\end{lemma}
It follows that $G$ induces a filter $K$, which is
$\Q$-generic over $V$. Moreover, it is clear that $V[K]=V[\vec\kappa, \vec H]$.

\subsection{$V_2$ is a non-trivial extension of $V_1$}
In this subsection we show that $\vec\eta \notin V_1$, and hence $V_1 \subsetneq V_2.$
Let $\sigma$ be an automorphism of $\Pp$ such that for each $p \in \Pp,$
\begin{itemize}
\item $\sigma$ fixes $\pi(p)$,

\item $\sigma$ fixes tail of $p,$

\item $\sigma$ permutes $\xi^p_i$ in the interval $(\nu^p_i, \nu^p_{i+1})$ which has no fixed points, for $i+1 < m^p$.
\end{itemize}
We then extend $\sigma$ to  some $\hat{\sigma}: V^{\Pp} \simeq V^{\Pp}$, and it then follows that $\hat{\sigma} \upharpoonright V^{\Q}$ is identity.
Now assume by the way of contradiction that $\vec\eta\in V_1$, and let $\dot{\vec\eta}$ be a $\Q$-name for it.
By density arguments, it is easily seen that the set
\[
D=\{p \in \Pp: \exists i< m^{p}: \xi^{p}_i \neq \sigma(\xi^{p}_i)   \}
\]
is dense in $\Pp$.
 Let $p \in G \cap D$ be a normal condition forcing
$\dot{\vec\eta} \in V_1$. Then:
\begin{enumerate}
\item As $\hat{\sigma} \upharpoonright V^{\Q}$ is identity, $\hat{\sigma}(\dot{\vec\eta})=\dot{\vec\eta}$,

\item As $p \in D$, $p \forces \hat{\sigma}(\dot{\vec\eta})\neq \dot{\vec\eta}$.
\end{enumerate}
By (1) and (2) we get a contradiction.

\subsection{Some canonization results-I}
 In this subsection we prove a few results  that will be used later. We concentrate here on forcing notions that appear
 in the stem of a condition.
We first introduce some auxiliary notations.
\begin{notation}
\begin{enumerate}
\item Assume $\vec\nu=\langle \nu_0, \cdots, \nu_{m-1}  \rangle$ is an increasing sequence from $X_0$. Then
\begin{enumerate}
\item $ \Sforce_{\vec\nu}=\prod_{i<m-1}\Col((\nu_i)^+,<\nu_{i+1})$,
\item $\R_{\vec\nu}=\Col((\nu_{m-1})^+, < \kappa)$.
\end{enumerate}
\item If  $p \in \Pp$, then  we set $ \Sforce_p=\Sforce_{\vec\nu_p}$, and $ \R_p=\R_{\vec\nu_p}$.

\end{enumerate}
\end{notation}
We have the following easy lemma.
\begin{lemma}
\label{someeasy}
Assume $p\in\Pp$, $t$ is legal for $p$ and $p\forces\dot a\subseteq\kappa$.
\begin{enumerate}
\item  There are a normal condition $q \leq^* p$ and a nice $\Sforce_{p}$-name $d$
such that $\vec f_{q}\upharpoonright m^p-1=\vec f_{p}\upharpoonright m^p-1$ and
$q\forces\dot a\cap\nu^q_{m^p-1}=d.$

\item If $q \leq^* p[t]$ and  $\vec f_{q}\upharpoonright m^{p[t]}-1=\vec f_{p[t]}\upharpoonright m^{p[t]}-1$, then there is normal condition $r \leq^* p$ such that
\begin{itemize}
\item $\stem(r)=\stem(q)\upharpoonright m^p$,
\item $\tail({r[t]}) \leq \tail(q)$,
\item $A^r\subseteq \big( A^p \cap (\xi^p_{m^p-1}+1)\big) \cup A^q$,
\item $B^r_\alpha \subseteq \big( B^p_\alpha \cap (\xi^p_{m^p-1}+1)\big) \cup B^q_\alpha$, for $0<\alpha < \kappa.$
\end{itemize}
Let $r:=r(p, t, q)$.
\end{enumerate}
\end{lemma}
\begin{proof}
clause (1) follows from Lemmas \ref{prikry property} and \ref{lem:factor}. Clause (2) is easy.
\end{proof}

\begin{definition}
\label{def:limit}
Assume $\langle p_i:i<\lambda\rangle$ is a $\leq^*$-decreasing sequence of conditions in $\Pp$ with the same stem. We define
$p:=\lim_{i<\lambda}p_i$ by $p=\stem(p)\concat{\tail(p)}$, where $\stem(p)=\stem(p_0)$, and $\tail(p)=(A, \langle B_\alpha: 0<\alpha < \kappa         \rangle, F_1, F_2)$, where
\begin{itemize}
\item $A=\bigcap_{i<\lambda}A^{p_i}$,
\item $B_\alpha=\bigcap_{i<\lambda}B_\alpha^{p_i}$,
\item $F_1(\nu)=\bigcup_{i<\lambda}F_1^{p_i}$,
\item $F_2(\nu, \xi)=\bigcup_{i<\lambda}F_2^{p_i}(\nu, \xi)$.
\end{itemize}
\end{definition}
Note that in general $p$ may not be a condition in $\Pp$.
The next easy fusion lemma gives conditions under which such a $p$ is a condition.
\begin{lemma}\label{sm:diagonal-fusion}
Assume $\tau$ is a limit ordinal,
$\langle p_i:i<\tau\rangle$ is a $\leq^*$-decreasing sequence of conditions in $\Pp$ with the same stem, and
$\theta:\tau\longrightarrow\kappa$. Assume further that:
\begin{enumerate}
\item For every $i<\tau$ and $0<\alpha<\kappa,$
$A^{p_{i+1}}\cap(\theta(i)+1)=A^{p_i}\cap(\theta(i)+1)$ and
$B^{p_{i+1}}_{\alpha}\cap(\theta(i)+1)
          =B^{p_{i}}_{\alpha}\cap(\theta(i)+1).$

\item For every $i<\tau$, $\mu\in A^{p_{i+1}}\cap\theta(i)$ and $\xi\in \bigcup_{0<\alpha < \tau}B^{p_{i+1}}_{\alpha}\setminus(\mu+1)$,
$F_1^{p_{i+1}}(\mu)=F_1^{p_i}(\mu)$ and
$ F_2^{p_{i+1}}(\mu,\xi)=F_2^{p_{i}}(\mu,\xi).$

\item For every $\omega\leq\delta<\kappa,~\bigl|\{i<\tau:\theta(i)\leq\delta\}\bigr|\leq|\delta|.$
\end{enumerate}
Then $p=\lim_{i<\tau}p_i \in \Pp$ and $p \leq^* p_i$ for all $i<\tau.$
\end{lemma}

 For the rest of this subsection, we assume the following:
 \begin{itemize}
\item $\nu<\kappa$ is inaccessible,
\item  $\Sforce$ is a forcing notion of  size at most
$\nu$,
\item $\R=\Col(\nu^+,<\kappa),$
\item $\R_\lambda=\Col(\nu^+,<\lambda),$ where $\nu^+ < \lambda < \kappa$ is inaccessible.
\end{itemize}

The poset $\Sforce\times\R$ is $\kappa$-c.c., hence by the results of Section 1, in the extension by $\Sforce\times\R$, the ideal generated by $U_\alpha$, namely
\[
 I_\alpha=
 \{D\subseteq\kappa:(\exists B\in U_\alpha)\ D\cap B=\varnothing\}
\]
is $(\nu^+, \kappa)$-good.

\begin{lemma}\label{sm:lift}
Assume $B\in U_\alpha$, $B \cap (\nu+1)=\emptyset$, and let
$T_\xi\in\Col(\nu^+,[\xi,\kappa))$ for $\xi\in B$.
Then
\begin{enumerate}
\item
$1_{\Sforce\times\R}\forces$``$C_T=\{\xi\in B:T_\xi\in\dot H\}\in I_\alpha^+$''.

\item Assume $(s, r) \in \Sforce\times\R$ and it forces ``$\dot D\subseteq C_T$ and
$\dot D\in I_\alpha^+$''. Then there are $(s', r') \leq (s, r)$,
$B'\subseteq B$ in $U_\alpha$, and conditions $T'_\xi \leq T_\xi$ in
$\Col(\nu^+,[\xi,\kappa))$ such that
$(s',r')\forces$``$\{\xi\in B':T'_\xi\in\dot H\}\subseteq\dot D$''.
\end{enumerate}
\end{lemma}
\begin{proof}
(1). Assume $(s, r) \in \Sforce\times\R$ and $B_0\in U_\alpha$. Let
$\sup\supp(r) < \xi\in B\cap B_0$, and set $(s', r')=(s, r \cup T_\xi)$. Then $(s', r')$ forces
 $\xi\in C_T\cap B_0$.

(2). Let $E$ be the set of all $\supp(r) < \xi \in B$ such that there is $(s_\xi, r_\xi) \leq (s, r \cup T_\xi)$ with $(s_\xi, r_\xi) \forces$``$\xi \in \dot D$''. It is easily seen that $E \in U_\alpha$, as otherwise $B \setminus E$ is forced to be disjoint to $\dot D.$ Let $B'=B \cap E$.
By shrinking $B'$, and using $|\Sforce|\leq \nu$ and the normality of $U_\alpha$, we may assume that:
\begin{itemize}
\item For some $s'$ and all $\xi \in B', s_\xi=s'$,
\item For some  $r' \in \R$ and all $\xi \in B', r_\xi \upharpoonright \xi=r'$.
\end{itemize}
Then for all $\xi \in B', r_\xi=r'\cup T'_\xi$, for some $T'_\xi \leq T_\xi$. Moreover,
clearly, $(s',r')\forces$``$\{\xi\in B':T'_\xi\in\dot H\}\subseteq\dot D$''.
\end{proof}

\begin{lemma}\label{sm:minima}
Assume $E$ is an equivalence relation on $\kappa$ and
$C\in I_\alpha^+$. Then either $C\cap[\beta]_E \in I_\alpha^+$ for some $\beta<\kappa$,
 or $C\cap\{\xi:\xi=\min[\xi]_E\} \in I_\alpha^+$.
\end{lemma}
\begin{proof}
Define $f: C \to \kappa$ by $f(\xi)=\min[\xi]_E \leq \xi.$ If $\{\xi \in C:\xi=\min[\xi]_E\} \in I_\alpha^+$, we are done.
Otherwise, $D=\{\xi \in C:\min[\xi]_E < \xi\} \in I_\alpha^+$. Thus $f$ is regressive on $D$, hence by Lemma \ref{lem:normalfacts}(3),  there is some
$\beta < \kappa$ such that $f^{-1}\{\beta\} \in I_\alpha^+$, i.e., $C\cap[\beta]_E \in I_\alpha^+$.
\end{proof}

\begin{lemma}\label{sm:stable}
Assume $r \in \R,$ $k<\omega$, and for each finite increasing sequence $\vec\lambda=(\lambda_1,\ldots,\lambda_k)$ of inaccessible
cardinals  above $\nu$,  $\dot z_{\vec\lambda}$ is an $\Sforce\times\R_{\lambda_1}$-name for an ordinal.
Then  there are $r'\leq r$,  a set $C\in U$, and an $\Sforce$-name $\dot f$ for an
ordinal function on $[C]^k$ such that for every $\vec\lambda\in[C]^k$,  $(1_{\Sforce}, r')\forces$``$\dot z_{\vec\lambda}=\dot f(\vec\lambda)$''.
Moreover, if the ordinal names are forced to be in a set $X$ of size $< \kappa,$ then
$\dot f(\vec\lambda)$ can be chosen so that it only depends on $\lambda_1$.
\end{lemma}
\begin{proof}
Let us assume that $\Sforce$ is trivial, as the general case can be proved similarly, using the fact that $|\Sforce| \leq \nu < \lambda_1$.
For each sequence $\vec\lambda$ as above with $\lambda_1 > \sup(\supp(r))$, let $r_{\vec\lambda} \leq r$ in $\R_{\lambda_1}$ and ordinal $\alpha_{\vec\lambda}$
be such that $r_{\vec\lambda}\forces$``$\dot z_{\vec\lambda}=\check \alpha_{\vec\lambda}$''.
As $|\R_{\lambda_1}| \leq \lambda_1 < \kappa,$ by Rowbottom's theorem, we can find $C_{\lambda_1} \in U$ and $r_{\lambda_1}$ such that for all $\vec\lambda$ as above starting with $\lambda_1$, if $(\lambda_2,\ldots,\lambda_k) \in [C_{\lambda_1}]^{k-1},$
then $r_{\vec\lambda}=r_{\lambda_1}$.
 Again as each $\lambda_1$ is inaccessible,  $\sup(\supp(r_{\lambda_1}) < \lambda_1,$ hence again by normality, we can find
 $D \in U$ and some fixed $r'$ such that for each $\lambda_1 \in D, r_{\lambda_1}=r'$. Now $r'$, $C=D \cap \Delta_{\lambda_1<\kappa}C_{\lambda_1}$
and the function $f$ defined as $f(\vec\lambda)=\alpha_{\vec\lambda}$ are as required. The moreover part follows by one more application of
Rowbottom's theorem.
\end{proof}

The canonization arguments below follow  previous section, but must account for the
interleaved collapses. In Section~1 the functions being canonized have
ground-model values, but here, their values are ordinal names, and equality
of their interpretations may depend on the collapse generic.
Consequently the corresponding equality kernel must, in general, be
represented by a name. But distinct  names
can interpret as the same set.
This forces us to define several names of equivalence relations, and
use them to get our canonization results, which provide uniform, coherent equality
tests without  fixing the collapse generic in advance.
we now proceed towards this.

\begin{lemma}\label{sm:compare}
Assume that for each $\xi<\lambda_1$,
$\dot z_{\xi,\vec\lambda}$ is an ordinal name as in
Lemma~\ref{sm:stable}. There are an
$\Sforce\times\R$-name $\dot E$ and $C\in U$ such that $\dot E$  is a name for an equivalence relation on
$\kappa$ and  whenever
$\vec\lambda\in[C]^k$ and $\xi,\zeta<\lambda_1$, then
$
 \bb{\xi\dot E\zeta}
 =\bb{\dot z_{\xi,\vec\lambda}=\dot z_{\zeta,\vec\lambda}},
$
where the Boolean values are taken in
$\RO(\Sforce\times\R_{\lambda_1})$.
\end{lemma}
\begin{proof}
For notational simplicity let us write
$\mathbb A=\RO(\mathbb S),$
$\mathbb B_\lambda=\RO(\mathbb S\times\mathbb R_\lambda),$
and
$\mathbb B=\RO(\mathbb S\times\mathbb R).$
Without loss of generality, $\mathbb A\subseteq\mathbb B_\lambda\subseteq\mathbb B.$
Let also  \(\widehat r\in\mathbb B\) denote the Boolean value corresponding to \((1_{\mathbb S},r)\), where \(r\in\mathbb R\).

For $\vec\lambda$ and $\xi,\zeta<\lambda_1$ as above, set
$
d_{\xi,\zeta,\vec\lambda}
=
\bigl\llbracket
\dot z_{\xi,\vec\lambda}
=
\dot z_{\zeta,\vec\lambda}
\bigr\rrbracket_{\mathbb B_{\lambda_1}}.
$
By Lemma \ref{sm:stable}, if $\xi,\zeta<\kappa$ and $r\in\mathbb R$, we can find $q \leq r$,
$H\in U$, and $c\in\mathbb A$ such that every member of \(H\) is an inaccessible cardinal above
$\max\{\nu^+,\xi,\zeta,\sup(\supp(q))\}$, and for every $\vec\lambda\in[H]^k,$
$
\widehat q\wedge d_{\xi,\zeta,\vec\lambda}
=
\widehat q\wedge c.
$
Let $A_{\xi,\zeta}$ be a maximal antichain of conditions $q\in\mathbb R$ for which there are $H$ and $c$ as above
and for each $a\in A_{\xi,\zeta}$, let $H_a\in U$ and $c_a\in\mathbb A$  witness it.
As
$\mathbb R$ is \(\kappa\)-c.c.,
$|A_{\xi,\zeta}|<\kappa,$ thus there is
$b_{\xi,\zeta}<\kappa$ such that
$
b_{\xi,\zeta}>\max\{\nu^+,\xi,\zeta\},$ and for every $a\in A_{\xi,\zeta},$
$\supp(a)\subseteq b_{\xi,\zeta}$.
Set
$
C_{\xi,\zeta}
=
\big(\bigcap_{a\in A_{\xi,\zeta}}H_a\big)
\setminus(b_{\xi,\zeta}+1) \in U,
$
and
$
e_{\xi,\zeta}
=
\bigvee_{a\in A_{\xi,\zeta}}
(\widehat a\wedge c_a).
$
It is easily seen that for every $\vec\lambda\in[C_{\xi,\zeta}]^k$,
$e_{\xi,\zeta}\in\mathbb B_{\lambda_1}$
and
$e_{\xi,\zeta}=d_{\xi,\zeta,\vec\lambda}.$
Now let $\dot E$ be such that  for all $\xi, \zeta < \lambda_1,$
$
\llbracket\xi\mathrel{\dot E}\zeta\rrbracket_{\mathbb B}
=e_{\xi,\zeta}.
$
Then easily,
\begin{center}
$1_{\mathbb B}\forces$``$\dot E$ is an equivalence relation on $\kappa$''.
\end{center}
 For each $\xi<\kappa$, let
$D_\xi=\Delta_{\zeta<\kappa}C_{\xi,\zeta} \in U,$
and set
$C=\Delta_{\xi<\kappa}D_\xi\in U.$
By shrinking $C$, we may assume it only contains inaccessible cardinals above $\nu$.
Then $C$ and $\dot E$ are as required.
\end{proof}

\subsection{Some canonization results-II}
In the previous subsection, our concentration was on forcings arising from the stem of conditions. In this subsection we also consider those arising from the tail of forcing conditions as well. For this subsection, we assume the following:
\begin{itemize}
\item $p$ is a normal condition of length $m \geq 1$,

\item $\ell\geq m$ and $r\geq1$,

\item $t_{\vec\nu, \vec\xi}=\langle (\nu_i, \xi_i): i< \ell             \rangle$ is so that
$t_{\vec\nu, \vec\xi}\upharpoonright [m, \ell)$ is legal for $p$,

\item $\vec\lambda\in[A^p]^r$ is increasing above $\nu_{\ell-1}$,

\item $\nu_\ell=\lambda_1 > \xi_{\ell-1}$,

\end{itemize}
\begin{definition}\label{sm:def:scalar}
Assume $p$ is a normal condition of length $m \geq 1$, $\ell\geq m$ and $r\geq1$.
\begin{enumerate}
\item $\Sigma^p_{\ell, r}$ is the collection of all  tuples $(\vec\nu, \vec\xi, \vec\lambda)$ satisfying the above hypotheses,

\item Assume $(\vec\nu, \vec\xi, \vec\lambda) \in \Sigma^p_{\ell, r}$. Then
\begin{itemize}
\item $\mathbb P(\vec\nu,\lambda_1)=\prod_{i<\ell}\Col(\nu_i^+,<\nu_{i+1})=\Sforce_{\vec\nu\concat\langle \lambda_1 \rangle}.$

\item $\Gamma^p_{\vec\nu, \vec\xi,\lambda_1}=\vec f_{p[t_{\vec\nu, \vec\xi}]}\in\mathbb P(\vec\nu,\lambda_1)$,

\end{itemize}
\item $\mathcal Z^p_{\ell, r}$ is the collection of all sequences $\langle \dot z_{\vec\nu, \vec\xi, \vec\lambda}: (\vec\nu, \vec\xi, \vec\lambda) \in \Sigma^p_{\ell, r}  \rangle$, where $\dot z_{\vec\nu, \vec\xi, \vec\lambda}$ is a
 $\mathbb P(\vec\nu,\lambda_1)$-name  for an ordinal. Note that since $\vec\xi$ is uniquely  determined by $\xi_{\ell-1}$,
 for $\eta=\xi_{\ell-1} < \lambda_1$ we may define $\dot z^{\vec\nu}_{\eta,\vec\lambda}= \dot z_{\vec\nu,\vec\xi,\vec\lambda}.$
 We set $\dot z^{\vec\nu}_{\eta,\vec\lambda}= \check0$, if there is no such $\vec\xi.$
 \end{enumerate}
 \end{definition}
The following lemma is a one-step canonization for ordinal names. It stabilizes their equality kernels and gives a constant-or-injective dichotomy along the last coordinate.
\begin{lemma}\label{sm:rows}
Assume $p$ is a normal condition, $m=m^p\geq1$, $\ell>m$, $r\geq1$ and
$
 \mathcal Z= \langle\dot z_{\vec\nu,\vec\xi,\vec\lambda}:
(\vec\nu,\vec\xi,\vec\lambda)\in\Sigma^p_{\ell,r}\rangle \in \mathcal Z^p_{\ell, r}.
$
Let $\vec\nu^-=\vec\nu\restr(\ell-1)$ and $\nu=\nu_{\ell-1}$.
Then there exist a normal condition $q\leq^*p$ with $\stem(q)=\stem(p)$ and
for each $ (\vec\nu^-,\vec\sigma, \langle\nu\rangle\concat\vec\lambda) \in\Sigma^q_{\ell-1,r+1},$
 names
$ \dot E_{\vec\nu},\dot E^{\lambda_1}_{\vec\nu},
 \dot\tau_{\vec\nu,\vec\sigma}$ and
$ \dot w_{\vec\nu,\vec\sigma,\vec\lambda}$
such that:
\begin{enumerate}
\item $\dot E_{\vec\nu}$ is an
$\Sforce_{\vec\nu}\times\R_{\vec\nu}$-name for an equivalence
relation on $\kappa$,  $\dot E^{\lambda_1}_{\vec\nu}$ is a
$\mathbb P(\vec\nu,\lambda_1)$-name for an equivalence relation
on $\lambda_1$, and for every  $\vec\lambda$ and
$\eta,\zeta<\lambda_1$,
\[
 \bb{\check\eta\mathrel{\dot E_{\vec\nu}}\check\zeta}_{\RO(\Sforce_{\vec\nu}\times\R_{\vec\nu})}=\bb{\check\eta\mathrel{\dot E^{\lambda_1}_{\vec\nu}}\check\zeta}_{\RO(\mathbb P(\vec\nu,\lambda_1))}=\bb{\dot z^{\vec\nu}_{\eta,\vec\lambda} =\dot z^{\vec\nu}_{\zeta,\vec\lambda}}_{\RO(\mathbb P(\vec\nu,\lambda_1))}.
 \]

In particular $\dot E_{\vec\nu}$ is independent of
$\vec\sigma$ and $\vec\lambda$.

\item  $\dot\tau_{\vec\nu,\vec\sigma}$ and
$\dot w_{\vec\nu,\vec\sigma,\vec\lambda}$ are
$\Sforce_{\vec\nu}$-names and
$
 1_{\Sforce_{\vec\nu}}\forces$``$
 \dot\tau_{\vec\nu,\vec\sigma}\in2$ and
$\dot w_{\vec\nu,\vec\sigma,\vec\lambda}$ is an ordinal''. Moreover, if
$(\vec\nu,\vec\sigma\concat\langle\eta\rangle,\vec\lambda)
\in\Sigma^q_{\ell,r}$, then
\[
 \Gamma^q_{\vec\nu,\vec\sigma\concat\langle\eta\rangle,\lambda_1}
 \forces_{\mathbb P(\vec\nu,\lambda_1)}
 \begin{cases}
 \dot z_{\vec\nu,\vec\sigma\concat\langle\eta\rangle,\vec\lambda}
       =\dot w_{\vec\nu,\vec\sigma,\vec\lambda},
                    &\text{if }\dot\tau_{\vec\nu,\vec\sigma}=0,\\
 (\forall\zeta<\check\eta)\,
       \neg(\zeta\mathrel{\dot E^{\lambda_1}_{\vec\nu}}\check\eta),
                    &\text{if }\dot\tau_{\vec\nu,\vec\sigma}=1.
 \end{cases}
\]

\item Assume $J\times H$ is
$\mathbb P(\vec\nu,\lambda_1)$-generic over $V$,
 $\vec\xi,\vec\zeta$ are distinct,
 $(\vec\nu,\vec\xi,\vec\lambda), (\vec\nu,\vec\zeta,\vec\lambda)\in\Sigma^q_{\ell,r}, $
 $ \Gamma^q_{\vec\nu,\vec\xi,\lambda_1},
 \Gamma^q_{\vec\nu,\vec\zeta,\lambda_1}\in J\times H, $ and
 $ \dot\tau_{\vec\nu,\vec\xi\restr(\ell-1)}[J]
   =\dot\tau_{\vec\nu,\vec\zeta\restr(\ell-1)}[J]=1.$
Then
 $
 \dot z_{\vec\nu,\vec\xi,\vec\lambda}[J\times H]
       \ne\dot z_{\vec\nu,\vec\zeta,\vec\lambda}[J\times H].
 $

\item The sequences $\langle \dot T_{\vec\nu^-,\vec\sigma, \langle\nu\rangle\concat\vec\lambda}: (\vec\nu^-,\vec\sigma, \langle\nu\rangle\concat\vec\lambda) \in \Sigma^q_{\ell-1, r+1}   \rangle$
    and
 $\langle \dot W_{\vec\nu^-,\vec\sigma,
                         \langle\nu\rangle\concat\vec\lambda}: (\vec\nu^-,\vec\sigma, \langle\nu\rangle\concat\vec\lambda) \in \Sigma^q_{\ell-1, r+1}   \rangle$
are in    $\mathcal Z^q_{\ell-1, r+1}$, where
$\dot T_{\vec\nu^-,\vec\sigma, \langle\nu\rangle\concat\vec\lambda}=\dot\tau_{\vec\nu,\vec\sigma},$
and $\dot W_{\vec\nu^-,\vec\sigma, \langle\nu\rangle\concat\vec\lambda}=\dot w_{\vec\nu,\vec\sigma,\vec\lambda}$.
\end{enumerate}
\end{lemma}
\begin{proof}
For each relevant $\vec\nu$, apply Lemma~\ref{sm:compare} to  $\vec\nu$ to find
$\dot E_{\vec\nu}$ and $C_{\vec\nu}\in U$ as there, and set
 $C_{\vec\nu}=\kappa$ otherwise. Let
$
 A^0=A^p\cap\{\lambda<\kappa:
       (\forall\vec\nu\in[\lambda]^\ell)\ \lambda\in C_{\vec\nu}\}\in U,
$
and let $p^0$ be obtained from $p$ by restricting $A^p$ to $A^0$.
For $\lambda_1\in A^0$,  $\vec\nu\in[\lambda_1]^\ell$, and
$\eta,\zeta<\lambda_1$, set
$
 b^{\vec\nu}_{\eta,\zeta}
 =\bb{\check\eta\mathrel{\dot E_{\vec\nu}}\check\zeta},$
and
$
 \dot E^{\lambda_1}_{\vec\nu}
 =\{(\check{(\eta,\zeta)},b^{\vec\nu}_{\eta,\zeta}):
                                     \eta,\zeta<\lambda_1\}.
$
By Lemma~\ref{sm:compare}, these Boolean values belong to
$\RO(\mathbb P(\vec\nu,\lambda_1))$, and clause~(1) holds below $p^0$.

For $\nu\in A^0$, let
$
 \mathcal I_\nu=\{(\vec a,\vec\sigma):
   (\exists\vec\lambda)\,
   (\vec a,\vec\sigma,\langle\nu\rangle\concat\vec\lambda)
                            \in\Sigma^{p^0}_{\ell-1,r+1}\}$
  and set $\mathcal I_\nu=\varnothing$ otherwise.
Clearly $|\mathcal I_\nu|\leq\nu$.  Define
$\mathcal I=j(\langle\mathcal I_\nu:\nu<\kappa\rangle)(\kappa)$.
For $x=(\vec a,\vec\sigma)\in\mathcal I$  set
\begin{itemize}
\item $\alpha(x)=\sigma_{\ell-2},$
\item $\mathbb S_x=\Sforce_{\vec a\concat\langle\kappa\rangle},$
\item $\widehat E_{\vec a}
   =j(\langle\dot E_{\vec\nu}\rangle)_{
                                  \vec a\concat\langle\kappa\rangle}.$
\end{itemize}
Then $\mathbb S_x\in V_{\kappa_1}^M$ and
$|\mathbb S_x|^M=\kappa$.
Let $\widehat z_{\vec a,\eta,\vec\lambda}$ denote the corresponding
 name from $j(\mathcal Z)$. Define the forcing notions
$\mathbb R_1=\Col(\kappa^+,<\kappa_1)^M,$
$\mathbb R_{2,\alpha} =\Col(\kappa^+,[\kappa_1,\kappa_{2,\alpha}))^{N_\alpha},$
and set
$\mathbb Q_\alpha=\mathbb R_1\times\mathbb R_{2,\alpha}.$ Clearly
$\mathbb Q_\alpha\simeq i_\alpha(\mathbb R_1)$, so we may identify them.
Now
\[
 1_{\mathbb S_x\times\mathbb Q_{\alpha(x)}}\forces^{N_{\alpha(x)}}
                          \dot u_x=\min[\check\kappa_1]_{i_{\alpha(x)}(\widehat E_{\vec a})},\leq\check\kappa_1.
\]
Since $K_1\times K_{2,\alpha(x)}$ is $\mathbb Q_{\alpha(x)}$-generic,
there are $a_x\in K_1$, $b_x\in K_{2,\alpha(x)}$, and
$\mathbb S_x$-names $\dot\tau_x,\dot\beta_x$ such that
\begin{itemize}
\item $a_x\leq j(F^{p^0}_1)(\kappa),$
\item  $b_x\leq k_{\alpha(x)}(F^{p^0}_2)(\kappa,\kappa_1),$
\item  $(1_{\mathbb S_x},(a_x,b_x))
       \forces^{N_{\alpha(x)}}_{\mathbb S_x\times\mathbb Q_{\alpha(x)}}
 \dot u_x=
 \begin{cases}
   \dot\beta_x,&\dot\tau_x=0,\\
   \check\kappa_1,&\dot\tau_x=1.
 \end{cases}
 $
 \end{itemize}
Then
$ \dot\tau_x,\dot\beta_x\in V_{\kappa_1}^{N_{\alpha(x)}} =V_{\kappa_1}^M,$
and
$ 1_{\mathbb S_x}\forces \dot\tau_x\in2\text{ and }\dot\beta_x<\check\delta_x$
By Lemma~\ref{sm:stable} in $M$ and the genericity of $K_1$, there are
$c_x\in K_1$, $C_x\in j(U)$, and $\mathbb S_x$-names
$\dot w_{x,\vec\lambda}$ such that
 \begin{itemize}
\item $c_x\leq a_x,$
\item $C_x\subseteq j(A^0)\setminus(\max\{\kappa,\delta_x\}+1),$
\item If $\vec\lambda\in[C_x]^r$, then $(1_{\mathbb S_x},c_x)
       \forces^M_{\mathbb S_x\times\mathbb R_1}
       \widehat z_{\vec a,\dot\beta_x,\vec\lambda}
                       =\dot w_{x,\vec\lambda}.$
 \end{itemize}
Set
\[
 r_1=j(F^{p^0}_1)(\kappa)\cup\bigcup_{x\in\mathcal I}c_x\in K_1
 \]
 and
 \[
 r_{2,\alpha}=k_\alpha(F^{p^0}_2)(\kappa,\kappa_1)
                \cup\bigcup_{\substack{x\in\mathcal I\\\alpha(x)=\alpha}}b_x
                                                    \in K_{2,\alpha}.
\]
Let $T_\alpha\in M$ be such that
$i_\alpha(T_\alpha)(\kappa_1)=r_{2,\alpha}$.
By \L o\'s's theorem and the $\kappa_1$-completeness of $j(U_\alpha)$,
there is $B^*_\alpha\in j(U_\alpha)$ such that, for
$\eta\in B^*_\alpha$ and $x\in\mathcal I$ with $\alpha(x)=\alpha$, the following conditions are satisfied:

 \begin{itemize}
\item $B^*_\alpha\subseteq j(B^{p^0}_\alpha)   \setminus(\max\{\kappa,\sup\supp(r_1)\}+1),$
\item $T_\alpha(\eta)\in\Col(\kappa^+,[\eta,\kappa_1))^M,$
\item $T_\alpha(\eta)\leq j(F^{p^0}_2)(\kappa,\eta),$
\item $(1_{\mathbb S_x},r_1\cup T_\alpha(\eta))
      \forces^M_{\mathbb S_x\times\mathbb R_1}
 \min[\check\eta]_{\widehat E_{\vec a}}
       =\begin{cases}
          \dot\beta_x,&\dot\tau_x=0,\\
          \check\eta,&\dot\tau_x=1.
        \end{cases}
        $
 \end{itemize}
Note that the sets $B^*_\alpha\subseteq j(X_\alpha)$ are pairwise disjoint.
For $\eta \in B^*_\alpha$ set $\widehat T(\eta)=T_\alpha(\eta)$.
and extend $\widehat T$ on the domain
$j(F^{p^0}_2)(\kappa,-)$ by $\widehat T(\eta)=j(F^{p^0}_2)(\kappa,\eta)$.
Clearly,
$
 \langle\dot\tau_x,C_x,
       \langle\dot w_{x,\vec\lambda}:\vec\lambda\in[C_x]^r\rangle:
                                      x\in\mathcal I\rangle, \langle T_\alpha,B^*_\alpha:0<\alpha<\kappa\rangle\in M.
$
By taking their $U$-representatives, we get $f,T$, names
$\dot\tau_{\nu,x},\dot w_{\nu,x,\vec\lambda}$, sets
$B_{\nu,\alpha}\in U_\alpha$, $C_{\nu,x}\in U$, and $A^1\in U$,
$A^1\subseteq A^0$, such that

 \begin{itemize}
\item $j(f)(\kappa)=r_1,$
\item $j(T)(\kappa,\cdot)=\widehat T,$
\item $k_\alpha(T)(\kappa,\kappa_1)=r_{2,\alpha},$
\item $f(\nu)\leq F^{p^0}_1(\nu),$
\item $B_{\nu,\alpha}\subseteq B^{p^0}_\alpha
          \setminus(\max\{\nu,\sup\supp(f(\nu))\}+1),$
\item $T(\nu,\eta)\leq F^{p^0}_2(\nu,\eta),$
\item $C_{\nu,x}\subseteq A^0\setminus(\nu+1)$, where
$\nu\in A^1,\ 0<\alpha<\nu,\ \eta\in B_{\nu,\alpha}$ and
$ x\in\mathcal I_\nu$.
 \end{itemize}
Note that if $\nu\in A^1$, $x=(\vec a,\vec\sigma)\in\mathcal I_\nu$,
$\eta\in B_{\nu,\alpha(x)}$, $\vec\lambda\in[C_{\nu,x}]^r$,
$\eta<\lambda_1$ and  $\vec\nu=\vec a\concat\langle\nu\rangle$, then
\[
 1_{\Sforce_{\vec\nu}}\forces
    \dot\tau_{\nu,x}\in2\text{ and }
    \dot w_{\nu,x,\vec\lambda}\in\Ord
\]
and
\[
 (1_{\Sforce_{\vec\nu}},f(\nu)\cup T(\nu,\eta))
       \forces_{\Sforce_{\vec\nu}\times\R_{\vec\nu}}
 \begin{cases}
  \dot z^{\vec\nu}_{\eta,\vec\lambda}
                   =\dot w_{\nu,x,\vec\lambda},
                                      &\dot\tau_{\nu,x}=0,\\
  \min[\check\eta]_{\dot E_{\vec\nu}}=\check\eta,
                                      &\dot\tau_{\nu,x}=1.
 \end{cases}
\]
Set
$ A={}A^1\cap\{\lambda<\kappa:
       (\forall\nu\in A^1\cap\lambda)(\forall x\in\mathcal I_\nu)
 \lambda\in C_{\nu,x}\},$
and
$ B_\alpha={}B^{p^0}_\alpha\cap\{\eta<\kappa:
       (\forall\nu\in A^1\cap\eta)
                    (\alpha<\nu\ \Longrightarrow\ \eta\in B_{\nu,\alpha})\}.
$
Then
$A\in U$ and $B_\alpha\in U_\alpha$.
Let $q$ be a condition such that  $\stem(q)=\stem(p)$ and
$\tail(q)=(A,\langle B_\alpha:0<\alpha<\kappa\rangle,f,T)$, with the functions
restricted to their required domains. Then
$q\leq^*p^0,$
$ j(F^q_1)(\kappa)=r_1\in K_1,$
and
$ k_\alpha(F^q_2)(\kappa,\kappa_1)=r_{2,\alpha}\in K_{2,\alpha}.$
By Lemma~\ref{ic:normal} we may assume that  $q$ is
normal. For every tuple in $\Sigma^q_{\ell,r}$,
$\eta\in B_{\nu,\alpha(x)}$ and $\vec\lambda\in[C_{\nu,x}]^r$
set
$
 \dot\tau_{\vec a\concat\langle\nu\rangle,\vec\sigma}
          =\dot\tau_{\nu,(\vec a,\vec\sigma)}$
          and
$ \dot w_{\vec a\concat\langle\nu\rangle,\vec\sigma,\vec\lambda}
          =\dot w_{\nu,(\vec a,\vec\sigma),\vec\lambda}.
$

Then $f(\nu)\cup T(\nu,\eta)\in\Col(\nu^+,<\lambda_1)$ and
\[
 \Gamma^q_{\vec\nu,\vec\sigma\concat\langle\eta\rangle,\lambda_1}
   =(\Gamma^q_{\vec a,\vec\sigma,\nu},f(\nu)\cup T(\nu,\eta))
   \leq(1_{\Sforce_{\vec\nu}},f(\nu)\cup T(\nu,\eta)).
\]
Clause (2) follows immediately.

For clause~(3), set $\eta=\xi_{\ell-1}$ and
$\theta=\zeta_{\ell-1}$. As $X_\alpha$'s are pairwise disjoint, if
$\vec\xi\ne\vec\zeta$, then $\eta\ne\theta$.
Assume $\eta<\theta$. By items~(1) and~(2),
\[
 \dot\tau_{\vec\nu,\vec\zeta\restr(\ell-1)}[J]=1
 \ \Longrightarrow\
 \neg\bigl(\eta\mathrel{\dot E^{\lambda_1}_{\vec\nu}[J\times H]}\theta\bigr)
 \ \Longrightarrow\
 \dot z_{\vec\nu,\vec\xi,\vec\lambda}[J\times H]
 \ne\dot z_{\vec\nu,\vec\zeta,\vec\lambda}[J\times H].
\]
Finally, clause~(4) follows as
$\mathbb P(\vec\nu^-,\nu)=\Sforce_{\vec\nu}$ and
$ 1_{\Sforce_{\vec\nu}}\forces
       \dot\tau_{\vec\nu,\vec\sigma}\in2\text{ and }
       \dot w_{\vec\nu,\vec\sigma,\vec\lambda}\in\Ord$
\end{proof}

\subsection{Minimality for sets of ordinals}
In this subsection, we show that $V_2$ is a minimal extension of $V_1$.
we first show that if  $a \subseteq \kappa$ is in $V_2 \setminus V_1$,
then $\vec\eta \in V_1[a].$

\begin{lemma}\label{sm:reading}
Assume $p\in\Pp$ and $p\forces\dot a\subseteq\kappa$.
Then there are a normal condition $q\leq^*p$ and a family
$
 \bigl\langle
 d_{\vec\nu_{q[t]},\,\vec\xi_{q[t]}\restr m^{q[t]}-1}:
 \varnothing\ne t\text{ is legal for }q
 \bigr\rangle$,
such that, for every  such $t$:
\begin{enumerate}
\item $d_{\vec\nu_{q[t]},\,\vec\xi_{q[t]}\restr m^{q[t]}-1}$ is a nice
$\Sforce_{q[t]}$-name for a subset of
$\nu^{q[t]}_{m^{q[t]}-1}$;
\item
$ q[t]\forces\dot a\cap\nu^{q[t]}_{m^{q[t]}-1}
       =d_{\vec\nu_{q[t]},\,\vec\xi_{q[t]}\restr m^{q[t]}-1}.$
\end{enumerate}
In particular, if $G$ is $\Pp$-generic over $V$, $q\in G$, and
$a=\dot a[G]$, then for all $n\geq m^p,$
\[
 a\cap\kappa_n
 =d_{\langle\kappa_0,\ldots,\kappa_n\rangle,
                   \vec\eta\restr n}
       [\langle H_i:i<n\rangle]
\]
\end{lemma}
\begin{proof}
For every normal condition $r \leq^* p$ such that $ \vec f_r=\vec f_{p}$,
let $\mathcal C(r)$ be the set of all pairs $(\vec\nu, \vec\xi)$, with $\lh(\vec\nu)=n+1, \lh(\vec\xi)=n$, such that for some legal
$t_{\vec\nu, \vec\xi}:=\langle(\nu_i,\xi_i):m^p\leq i < n\rangle$ for $r$, we have
$(\vec\nu \upharpoonright n, \vec\xi)=(\vec\nu_{r[t_{\vec\nu, \vec\xi}]}, \vec\xi_{r[t_{\vec\nu, \vec\xi}]})$
and $\xi_{n-1} < \nu_n \in A^r$.

For  such an $r$ and $(\vec\nu, \vec\xi) \in \mathcal C(r)$, we recursively define a $\leq^*$-decreasing sequence
$\langle r^\delta:\delta\leq\kappa\rangle$ of conditions as follows
Set $r^0=r$.
Assume now that $r^\xi$ is defined. If
${\nu_n}<\xi\in B^{r^\xi}_{\xi_{n-1}}$, set
$t_\xi=t_{\vec\nu, \vec\xi}\concat\langle({\nu_n},\xi)\rangle$. Then by Lemma \ref{someeasy}(1), there are
$q^\xi \leq^* r^\xi[t_\xi]$ and an $\Sforce_{\vec\nu}$-name $d_\xi$ such that
$ q^\xi\forces\dot a\cap{\nu_n}=d_\xi$ and
$\vec f_{q^\xi}\restr n=\vec f_{r[t]} \upharpoonright n \in \Sforce_{\vec\nu}.$
Let $r^{\xi+1}:=r(r^\xi, q^\xi, t_\xi)$, as in Lemma \ref{someeasy}(2). Set also $f_\xi=f^{q^\xi}_{n}$.
Otherwise set $r^{\xi+1}=r^\xi$.
At limit stages, set $r^\xi=\lim_{\zeta<\xi}r^\zeta$. By our construction and Lemma \ref{sm:diagonal-fusion},  $r^\xi$ is a condition.
Let $r_*=r^\kappa$.

 Then for every ${\nu_n}<\xi\in B^{r_*}_{\xi_{n-1}}$, we have $\vec f_{r_*[t_\xi]}\restr n=\vec f_{r[t]} \upharpoonright n$, and
$\tail(r_*[t_\xi])\leq\tail(q^\xi).$
Since $|\Sforce_{\nu_n}|\leq\nu_n<\kappa$, there are
$C\subseteq B^{r_*}_{\xi_{n-1}} \setminus \nu_n+1$ in $U_{\xi_{n-1}}$ and a nice $\Sforce_{\vec\nu}$-name
$d_{\vec\nu,\vec\xi}$ such that $d_\xi=d_{\vec\nu,\vec\xi}$ for all $\xi \in C.$
By further shrinking $C$, we may assume that for some $g \in \Col(\nu_n^+, < \kappa)$ and  all $\xi \in C$,
$f_\xi \upharpoonright \xi=g.$
Let $\bar r:= \bar r_{\vec\nu, \vec\xi} \leq^* r^*$ be a normal condition such that
\begin{itemize}
\item $\vec f_{\bar r}=\vec f_{r^*}$,
\item $B^{\bar r}_{\xi_{n-1}}=(B^{r_*}_{\xi_{n-1}}\cap(\nu_n+1))\cup C,$
\item $B^{\bar r}_\beta=(B^{r_*}_\beta\cap(\nu_n+1))
             \cup(B^{r_*}_\beta\setminus(\gamma+1))$, where $\gamma=\max\{\nu_n,\sup(\supp(g))\}$,
\item $F_1^{\bar r}(\nu_n)=g$,
\item $F_2^{\bar r}(\nu_n, \xi)=f_\xi \upharpoonright [\xi, \kappa)$,
\end{itemize}
Now
if $\nu_n<\xi\in B^{\bar r}_{\xi_{n-1}}$, then
$
 \bar r[t_\xi]\leq q_\xi,$ and hence for all such $\xi$ we have
\[
 \bar r[t_{\vec\nu, \vec\xi}\concat\langle(\nu_n,\xi)\rangle]
       \forces\dot a\cap\nu=d_{\vec\nu, \vec\xi}.
\]
Now, given $r$ and $\nu\in A^r$, the set $ \mathcal C_\nu(r)=\{(\vec\nu, \vec\xi) \in  \mathcal C(r):\nu_n=\nu               \}$
has size at most $\nu$, hence we can repeat the above construction and obtain a condition $\bar r(\nu)$ which is $\leq^* \bar r_{\vec\nu, \vec\xi}$,
for all $(\vec\nu, \vec\xi) \in  \mathcal C_\nu(r)$

Again, by induction on $\nu < \kappa$, we can form a $\leq^*$-decreasing sequence of normal conditions $\langle r(\nu): \nu \leq \kappa     \rangle$
such that $q=r(\kappa) \leq^* p$ is a well-defined condition, and for every $(\vec\nu, \vec\xi) \in \mathcal C(q_0)$ and
every $\nu_n<\xi\in B^{q_0}_{\xi_{n-1}}$, we have
\[
 q[t_{\vec\nu, \vec\xi}\concat\langle(\nu_n,\xi)\rangle]
       \forces\dot a\cap\nu_n=d_{\vec\nu, \vec\xi}.
\]
It is clear that $q$ and the sequence $
 \bigl\langle
 d_{\vec\nu_{q[t]},\,\vec\xi_{q[t]}\restr m^{q[t]}-1}:
 \varnothing\ne t\text{ is legal for }q
 \bigr\rangle$ are as required. The last part of the lemma follows immediately.
\end{proof}
We have the following immediate corollary.
\begin{corollary}\label{sm:absorption}
Every bounded subset of $\kappa$ in $V[G]$ belongs to $V_1$, and
every subset of $\kappa$ in $V[G]$ belongs to $V_2$.
\end{corollary}
\begin{definition}
Assume  $p$ is a normal condition, $m=m^p\geq1$, $\ell\geq m$, $r\geq1$, and
$
 \mathcal Z=
 \langle\dot z_{\vec\nu,\vec\xi,\vec\lambda}:
       (\vec\nu,\vec\xi,\vec\lambda)\in\Sigma^p_{\ell,r}\rangle
       \in\mathcal Z^p_{\ell,r}.
$
Assume $G$ is
$\mathbb P(\vec\nu,\lambda_1)$-generic  over $V$. Then
\begin{itemize}
\item $D^p_{\vec\nu,\vec\lambda}(G)=
 \left\{\vec\xi:
       (\vec\nu,\vec\xi,\vec\lambda)\in\Sigma^p_{\ell, r},
       \Gamma^p_{\vec\nu,\vec\xi,\lambda_1}\in G \right\},
$
\item $G_p=G\restr\Sforce_p.$
\end{itemize}
\end{definition}

\begin{lemma}\label{sm:canon}
Assume $p$ is a normal condition, $m=m^p\geq1$, $\ell\geq m$, $r\geq1$, and
$
 \mathcal Z=
 \langle\dot z_{\vec\nu,\vec\xi,\vec\lambda}:
       (\vec\nu,\vec\xi,\vec\lambda)\in\Sigma^p_{\ell,r}\rangle
       \in\mathcal Z^p_{\ell,r}.
$
Then there exist a normal condition $q\leq^*p$ with $ \vec f_q\restr(m-1)=\vec f_p\restr(m-1)$, an $\Sforce_p$-name $\dot k$ for an integer
in $\{m,\ldots,\ell\}$, and
$\mathbb P(\vec\nu,\lambda_1)$-names
$\langle \dot h_{\vec\nu,\vec\lambda}: \exists \vec\xi,~(\vec\nu, \vec\xi, \vec\lambda) \in \Sigma^q_{\ell,r} \rangle$
such that:
\begin{enumerate}
\item If $G$ is
$\mathbb P(\vec\nu,\lambda_1)$-generic  over $V$,
 $k=\dot k[G_p]$ and
$D=D^q_{\vec\nu,\vec\lambda}(G)$, then $ \dot h_{\vec\nu,\vec\lambda}[G]$ is an injective function with
$ \dom\bigl(\dot h_{\vec\nu,\vec\lambda}[G]\bigr)
  =\{\vec\xi\restr k:\vec\xi\in D\}$
and for all $\vec\xi\in D,$
$ \dot z_{\vec\nu,\vec\xi,\vec\lambda}[G]
   =\dot h_{\vec\nu,\vec\lambda}[G](\vec\xi\restr k).
$
In particular, for all $\vec\xi,\vec\zeta\in D$,
\[
 \dot z_{\vec\nu,\vec\xi,\vec\lambda}[G]
   =\dot z_{\vec\nu,\vec\zeta,\vec\lambda}[G]
 \quad\Longleftrightarrow\quad
 \vec\xi\restr k=\vec\zeta\restr k.
\]
\item If $X$ is a finite set of ordinals such that
$
 \Gamma^p_{\vec\nu,\vec\xi,\lambda_1}
       \forces_{\mathbb P(\vec\nu,\lambda_1)}
       \dot z_{\vec\nu,\vec\xi,\vec\lambda}\in\check X$
 for all $\vec\nu,\vec\xi,\vec\lambda)
                   \in\Sigma^p_{\ell,r}$,
then we can chose $q$ as in (1) and an
$\Sforce_p$-name $\dot c$ such that  $\dot k=\check m$,
$
 1_{\Sforce_p}\forces\dot c\in\check X,
$
and, for every $(\vec\nu,\vec\xi,\vec\lambda)\in\Sigma^q_{\ell,r}$,
$
 \Gamma^q_{\vec\nu,\vec\xi,\lambda_1}
       \forces_{\mathbb P(\vec\nu,\lambda_1)}
       \dot z_{\vec\nu,\vec\xi,\vec\lambda}=\dot c.
$
\end{enumerate}
\end{lemma}
\begin{proof}
We prove (1) and (2) simultaneously, by induction on $\ell \geq m$ and
for all $r\geq1$.

First assume that $\ell=m$. Then $\{\vec\xi: (\vec\nu,\vec\xi,\vec\lambda) \in \Sigma^p_{\ell, r}   \}=\{\vec\xi_p \}$,
hence, for clause (1),  it suffices to take $q=p$ and
$\dot k=\check m$. Clause (2) follows easily from Lemma~\ref{sm:stable}. Indeed, by applying the lemma
to  $\Sforce=\Sforce_p,$
$ \R=\R_p=\Col(\rho^+,<\kappa)$ and
$ r=f^p_{m-1},$
we can find
$r'\leq r$, $C\in U$, and
$\Sforce_p$-names $\dot b_\lambda$ such that for every $\vec\lambda\in[C\cap A^p]^r,$
$ (1_{\Sforce_p},r_1)\forces
 \dot z_{\vec\nu_p,\vec\xi_p,\vec\lambda}  =\dot b_{\lambda_1}$.
Each $\dot b_\lambda$ is forced to be in $X$, hence by further shrinking $C$, we can find
$C'\in U$, $C'\subseteq C\cap A^p$,
and an $\Sforce_p$-name $\dot c$ such that for all $\lambda\in C'$,
$ 1_{\Sforce_p}\forces
       \dot b_\lambda=\dot c$.
Let $q \leq^* p$ be a normal condition which is obtained by
replacing $f^p_{m-1}$ by $r'$ and restricting $A^p$ to
$C'\setminus(\max\{\xi^p_{m-1},\sup\supp(r_1)\}+1)$.
Then this $q$ and $\dot c$ are as required.

Now assume $\ell>m$ and  items (1) and (2) hold for $\ell-1$.
We first prove clause~(2).
By Lemma~\ref{sm:rows}, there is  a normal condition $p'\leq^*p$ and the
names $ \dot E_{\vec\nu},\dot E^{\lambda_1}_{\vec\nu},
 \dot\tau_{\vec\nu,\vec\sigma}$ and
$ \dot w_{\vec\nu,\vec\sigma,\vec\lambda}$ as there.
For each $\vec\nu$, by clause~(1), $\Sforce_{\vec\nu}\times\R_{\vec\nu}$ forces ``$ |\kappa/\dot E_{\vec\nu}|\leq |X|$'',
hence the set
\[
 \dot M_{\vec\nu}=
       \{\eta<\kappa:\eta=\min[\eta]_{\dot E_{\vec\nu}}\}
\]
is forced to be finite.
Now assume that
$(\vec\nu^-,\vec\sigma,\langle\nu\rangle\concat\vec\lambda)
 \in\Sigma^a_{\ell-1,r+1}$. We claim that
$\Gamma^{p'}_{\vec\nu^-,\vec\sigma,\nu_{\ell-1}}
       \forces_{\Sforce_{\vec\nu}}
       \dot\tau_{\vec\nu,\vec\sigma}=0$.
 Assume not, and let  $J\times H$ be $\Sforce_{\vec\nu}\times\R_{\vec\nu}$-generic
 over $V$ such that $(\Gamma^{p'}_{\vec\nu^-,\vec\sigma,\nu_{\ell-1}}, F_1^{p'}(\nu_{\ell-1})) \in J \times H$
 and  $\dot\tau_{\vec\nu,\vec\sigma}[J]=1$. Set
       $B=B^{p'}_{\sigma_{\ell-2}}\setminus(\nu_{\ell-1}+1)$,
$f=F_1^{p'}(\nu_{\ell-1})$, and $T_\eta=F_2^{p'}(\nu_{\ell-1},\eta)$ for $\eta\in B$.
 By Lemma~\ref{sm:lift}(1),
$
 C_T=\{\eta\in B:T_\eta\in H\}\in I_\alpha^+.
$
Let $\eta\in C_T$ with $\lambda_1>\eta$.
Then $f\cup T_\eta \in \Col(\nu^+,<\lambda_1)$, and
\[
 \Gamma^{p'}_{\vec\nu,\vec\sigma\concat\langle\eta\rangle,\lambda_1}
       =(\Gamma^{p'}_{\vec\nu^-,\vec\sigma,\nu_{\ell-1}},f\cup T_\eta)
       \in (J\times H)\restr\mathbb P(\vec\nu,\lambda_1).
\]
Thus by Lemma~\ref{sm:rows},
$\eta\in\dot M_{\vec\nu}[J\times H]$.
Thus $C_T \subseteq M_{\vec\nu}[J\times H]$ is finite, which contradicts
$C_T \in I_\alpha^+$-positivity.

Now set
\[
 \dot W^c_{\vec\nu^-,\vec\sigma,
                         \langle\nu\rangle\concat\vec\lambda}
 =\begin{cases}
   \dot w_{\vec\nu,\vec\sigma,\vec\lambda},
       &\dot w_{\vec\nu,\vec\sigma,\vec\lambda}\in X,\\
   \check0,&\text{otherwise}.
 \end{cases}
\]
By the induction hypothesis,  there are a normal condition $q\leq^*p'$ and an
$\Sforce_p$-name $\dot c$. It is easily seen that $q$ and $\dot c$ are as required.

Now we prove clause~(1). again using Lemma~\ref{sm:rows}, we get   a normal condition $p'\leq^*p$ and names for equivalence relations as there.
By the induction hypothesis, applied to $\langle \dot T_{\vec\nu^-,\vec\sigma, \langle\nu\rangle\concat\vec\lambda}: (\vec\nu^-,\vec\sigma, \langle\nu\rangle\concat\vec\lambda) \in \Sigma^q_{\ell-1, r+1}   \rangle$, see Lemma~\ref{sm:rows}(4),
we get  a normal condition
$p''\leq^*p'$ and an $\Sforce_p$-name $\dot\tau\in2$. Using the induction hypothesis, applied to
 $\dot W$ as defined above, we get a further extension
 $q\leq^*p''$ and an $\Sforce_p$-name $\dot k'$, which is forced to be in
$\{ m, \cdots, \ell-1\}$.
Define the $\Sforce_p$-name $\dot k$ by $\dot k=\dot k'$ if $\tau=0$ and $\dot k=\check \ell$ if $\tau=1.$
We show that $q$ and $\dot k$ are as required.

Let  $(\vec\nu,\vec\lambda), k, G$ and $D$ be as clause (1), and set
$J=G\restr\Sforce_{\vec\nu}$.
If $\vec\xi\in D$ and
$\vec\sigma=\vec\xi\restr(\ell-1)$, then clearly $\Gamma^q_{\vec\nu,\vec\xi,\lambda_1}
   \restr(\ell-1) =\Gamma^q_{\vec\nu^-,\vec\sigma,\nu},$
   and $\mathbb P(\vec\nu^-,\nu)=\Sforce_{\vec\nu}$, hence
 $\dot\tau_{\vec\nu,\vec\sigma}[J]=\dot\tau[G_p]$.
If this common value is $0$, then we apply Lemma~\ref{sm:rows}(2), and the  induction
hypothesis for $\dot W$, and if this common value is $1$, then we apply Lemma~\ref{sm:rows}(3)
to get the desired result.
\end{proof}

Let $G$ be a $\Pp$-generic filter over $V$, $\pi: \PP \to \Q$ be the projection map and let $K$ be the
$\Q$-generic filter over $V$ induced from $\pi$. Let also $ \vec\nu_n=(\kappa_0,\ldots,\kappa_n)$
 and
$ J_n=\langle H_i:i<n\rangle.$
Recall from Lemma~\ref{sm:reading} that  $ a\cap\kappa_n=d_{v_n,\vec\eta\restr n}[J_n]$.

\begin{lemma}\label{sm:tables}
Assume $p\in\Pp$ is a normal condition, $m=m^p\geq1$,
$p\forces\dot a\subseteq\kappa$ and
$
 \bigl\langle
 d_{\vec\nu_{p[t]},\,\vec\xi_{p[t]}\restr m^{p[t]}-1}:
 \varnothing\ne t\text{ is legal for }p
 \bigr\rangle$
 satisfies the conclusions of
Lemma~\ref{sm:reading} for $p$ and $\dot a$.
Then there exist a normal condition $q\leq^*p$ with $ \vec f_q\restr(m-1)=\vec f_p\restr(m-1)$, a sequence
$\langle\dot k_n:n\geq m\rangle$ of $\Sforce_p$-names
with $1_{\Sforce_p}\forces\dot k_n\in\{m,\ldots,n\}$, and
a $\Q$-name $\dot{\mathcal T}$ such that:
\begin{enumerate}
\item If $\pi(q)\in K$, then letting $k_n=\dot k_n[J_{m-1}],$ we have
$
 \dot{\mathcal T}[K]=\langle k_n,D_n,L_n:n\geq m\rangle\in V_1,
$
where, for every $n\geq m$,
\begin{itemize}
\item $E_n=D^q_{\vec\nu_{n-1},\langle\kappa_n\rangle}(J_n),$
\item $D_n=\{\vec\xi\restr k_n:\vec\xi\in E_n\},$
\item $L_n:D_n\longrightarrow\pow(\kappa_n)^{V_1}$
is injective.
 \end{itemize}
In particular, $m\leq k_n\leq n$ and $D_n\subseteq{}^{k_n}\kappa$.
\item If  $q\in G$ and  $a=\dot a[G]$, then for every $n \geq m,$ $\vec\eta\restr k_n\in D_n$
and
$ L_n(\vec\eta\restr k_n)=a\cap\kappa_n.$
\end{enumerate}
\end{lemma}
\begin{proof}
For $n\geq m$ and for
$\vec v=\langle\nu_0,\ldots,\nu_n\rangle$ extending $\vec\nu_p$,
let
$
 \langle e^{\vec v}_\zeta:\zeta<\theta_{\vec\nu}\rangle
$
be an enumeration of
 all nice $\Sforce_{\vec v}$-names for subsets of $\nu_n$. Note that
 $0<\theta_{\vec v}\leq 2^{\nu_n}<\kappa$.
For $(\vec\nu,\vec\xi,\langle\lambda\rangle)\in\Sigma^p_{n,1}$,
set ${\vec v}=\vec\nu\concat\langle\lambda\rangle$ and let
$\dot c_{\vec\nu,\vec\xi,\langle\lambda\rangle}$ be the
$\mathbb P(\vec\nu,\lambda)=\Sforce_{\vec v}$-name which is forced to be
$
 \dot c_{\vec\nu,\vec\xi,\langle\lambda\rangle}
   =\min\{\zeta<\theta_{\vec v}:e^{\vec v}_\zeta=d_{{\vec v},\vec\xi}\}.
$
 The resulting sequence $\langle \dot c_{\vec\nu,\vec\xi,\langle\lambda\rangle}: (\vec\nu,\vec\xi,\langle\lambda\rangle)\in\Sigma^p_{n,1} \rangle \in \mathcal Z^p_{n,1}$. Moreover, for any two indices with the same
$(\vec\nu,\lambda)$, it is forced that
\[
 \dot c_{\vec\nu,\vec\xi,\langle\lambda\rangle}
     =\dot c_{\vec\nu,\vec\zeta,\langle\lambda\rangle}
 \quad\Longleftrightarrow\quad
 d_{\vec v,\vec\xi}=d_{\vec v,\vec\zeta}.
\]
We define a $\leq^*$-decreasing sequence of normal conditions with the same $\vec f_p\restr(m-1)$ and a sequence $\langle \dot k_n: n \geq m \rangle$
with  $1_{\Sforce_p}\forces\dot k_n\in\{m,\ldots,n\}$ as follows.
Set $p_m=p$. Now assume $p_n$ and $\dot k_{n-1}$ are defined.
By Lemma~\ref{sm:canon} applied to the above sequence of names restricted to
$\Sigma^{p_n}_{n,1}$, we get a normal condition $p_{n+1}\leq^*p_n$,
an $\Sforce_p$-name $\dot k_n$, and names
$\dot h^n_{\vec\nu,\langle\lambda\rangle}$ as there.
Let $q=\lim_{m \leq n<\omega}p_n$, which we may assume is normal.
Set
\begin{itemize}
\item $k_n=\dot k_n[J_{m-1}],$
\item $E_n=D^q_{\vec\nu_n\restr n,\,\langle\kappa_n\rangle}(J_n),$
\item $D_n=\{\vec\xi\restr k_n:\vec\xi\in E_n\},$
\item $h_n=\dot h^n_{\vec\nu_n\restr n,\,\langle\kappa_n\rangle}[J_n]
                                                   \restr D_n,$
\item $L_n(t)=e^{\vec\nu_n}_{h_n(t)}[J_n]$, where $t\in D_n$.
 \end{itemize}
It is not difficult to show that these are as required.
\end{proof}

\subsection{Completing the proof of Theorem \ref{thm:main2}}
We now complete the proof of Theorem \ref{thm:main2}.
Let $V_1=V[\vec\kappa, \vec H]$ and $V_2=V_1[\vec\eta]$ be as above. By Subsection 2.2, $V_1\subsetneq V_2$, and
by the above results, $V_2$ is a minimal extension of $V_1$ with respect to subsets of
$\kappa$.
Thus by  Lemmas \ref{lem:subsetreduction} and \ref{lem:chain}  $V_2$ is minimal over $V_1$.

\end{document}